\documentclass{gtart}

\def\ifplaintex{\expandafter\ifx\csname documentclass\endcsname\relax}

\def\gtp{{\mathsurround=0pt\it $\cal G\mskip-2mu$eometry \&\ 
$\cal T\!\!$opology $\cal P\!$ublications}}  

\def\recd{{\small Received:\qua\receiveddate\ifx\reviseddate\relax
\else\qquad Revised:\qua\reviseddate\fi\par}} 

\def\volumenumber#1{\def\thevolumenumber{#1}}
\def\volumeyear#1{\def\thevolumeyear{#1}}
\def\papernumber#1{\def\thepapernumber{#1}}
\def\pagenumbers#1#2{\def\startpage{#1}\def\finishpage{#2}}
\def\published#1{\def\publishdate{#1}}

\def\received#1{\def\receiveddate{#1}}

\def\asciiaddress#1{\def\theasciiaddress{#1}}

\long\def\asciiabstract#1{\long\def\theasciiabstract{#1}}

\let\\\par\let\thevolumenumber\relax
\let\thepapernumber\relax\let\thevolumeyear\relax\let\startpage\relax
\let\finishpage\relax\let\publishdate\relax\let\receiveddate\relax
\let\reviseddate\relax\let\theasciititle\relax
\let\theasciiauthors\relax\let\theasciiaddress\relax
\let\theasciiabstract\relax

\let\theasciiemail\relax

\ifplaintex
\font\logobig=cmssbx10 scaled 3836
\font\logomed=cmssbx10 scaled 2557
\else
\font\logobig=cmssbx10 scaled 4200
\font\logomed=cmssbx10 scaled 2800
\fi

\long\def\makeagttitle{   
\count0=\startpage
\agt\hfill      
\hbox to 45truept{\vbox to 0pt{\vglue -13truept{\logomed A\kern -.37em{\logobig 
T}\kern -.38em G}\vss}\hss}
\break
{\small Volume \thevolumenumber\ (\thevolumeyear)
\startpage--\finishpage\nl
Published: \publishdate}

\vglue .25truein

{\parskip=0pt\leftskip 0pt plus
1fil\def\\{\par\smallskip}{\Large\bf\thetitle}\par\medskip} \vglue
0.05truein

{\parskip=0pt\leftskip 0pt plus 1fil\def\\{\par}{\sc\theauthors}
\par\medskip}%
 
\vglue 0.03truein 

{\small\leftskip 25truept\rightskip 25truept{\bf Abstract}\stdspace\theabstract

{\bf AMS Classification}\stdspace\theprimaryclass
\ifx\thesecondaryclass\relax\else; \thesecondaryclass\fi\par
{\bf Keywords}\stdspace \thekeywords\par}\vglue 7truept

}   

\ifplaintex
\font\phead=cmsl9 scaled 950
\font\pnum=cmbx10 scaled 913
\font\pfoot=cmsl9 scaled 950
\headline{\vbox to 0pt{\vskip -4.5mm\line{\small\phead\ifnum
\count0=\startpage ISSN 1472-2739 (on-line) 1472-2747 (printed)
\hfill {\pnum\folio}\else\ifodd\count0\def\\{ }%
\ifx\theshorttitle\relax\thetitle\else\theshorttitle\fi\hfill{\pnum\folio}
\else\def\\{ and }{\pnum\folio}\hfill\ifx\theshortauthors\relax\theauthors
\else\theshortauthors\fi\fi\fi}\vss}}
\footline{\vbox to 0pt{\vglue 0mm\line{\small\pfoot\ifnum\count0=\startpage
\copyright\ \gtp\hfill\else
\agt, Volume \thevolumenumber\ (\thevolumeyear)\hfill\fi}\vss}}
\else
\font=cmsl9 scaled 1050
\font\lnum=cmbx10 
\font=cmsl9 scaled 1050
\makeatletter
\def\@oddhead{{\small\lhead\ifnum\count0=\startpage ISSN 1472-2739 
(on-line) 1472-2747 (printed)\hfill {\lnum\number\count0}\else\ifodd\count0
\def\\{ }\ifx\theshorttitle\relax \thetitle \else\theshorttitle\fi\hfill
{\lnum\number\count0}\else\def\\{ and }{\lnum\number\count0}
\hfill\ifx\theshortauthors\relax 
\theauthors\else\theshortauthors\fi\fi\fi}}\def\@evenhead{\@oddhead}
\def\@oddfoot{\small\lfoot\ifnum\count0=\startpage\copyright\ \gtp\hfill\else
\agt, Volume \thevolumenumber\ (\thevolumeyear)\hfill\fi}
\def\@evenfoot{\@oddfoot}
\makeatother
\fi
\let\maketitlepage\makeagttitle
\let\makeshorttitle\maketitlepage
\let\maketitle\maketitlepage

\newwrite\gtoutfile
\long\gdef\makeheadfile{  
{\def\\{, }\def\s{ }
\immediate\openout\gtoutfile head.xxx
\immediate\write\gtoutfile{To: math@arxiv.org}
\immediate\write\gtoutfile{Subject: put OR rep NNNNN:ppppp}
\immediate\write\gtoutfile{--text follows this line--}
\immediate\write\gtoutfile{Proxy-for: \ifx\theasciiauthors\relax
\theauthors\else\theasciiauthors\fi\s<\ifx\theasciiemail\relax\theemail\else\theasciiemail\fi>}
\immediate\write\gtoutfile{\noexpand\\}
\immediate\write\gtoutfile{Authors: \ifx\theasciiauthors\relax
\theauthors\else\theasciiauthors\fi}
{\def\\{ }\immediate\write\gtoutfile{Title: \ifx\theasciititle\relax
\thetitle\else\theasciititle\fi}}
\immediate\write\gtoutfile{Subj-class: GT or SG, GR etc}
\immediate\write\gtoutfile{MSC-class: \theprimaryclass\ifx\thesecondaryclass\relax\else, \thesecondaryclass\fi}
\immediate\write\gtoutfile{Journal-ref: Algebraic and Geometric Topology \thevolumenumber\s
(\thevolumeyear) \startpage-\finishpage}
\immediate\write\gtoutfile{Comments: Published in Algebraic and
Geometric Topology at}
\immediate\write\gtoutfile{\s\s\s  http://www.maths.warwick.ac.uk/agt/AGTVol\thevolumenumber/agt-\thevolumenumber-\thepapernumber.abs.html}
\immediate\write\gtoutfile{\noexpand\\}
\immediate\write\gtoutfile{}
\ifx\theasciiabstract\relax
\immediate\write\gtoutfile{\theabstract}\else
\immediate\write\gtoutfile{\theasciiabstract}\fi
\immediate\write\gtoutfile{}
\immediate\write\gtoutfile{\noexpand\\}
\immediate\write\gtoutfile{}
\immediate\closeout\gtoutfile}}  

\def\maketitlepage{\makeagttitle\makeheadfile}
\let\makeshorttitle\maketitlepage
\let\maketitle\maketitlepage

\def\ifplaintex{\expandafter\ifx\csname documentclass\endcsname\relax}

\def\gtp{{\mathsurround=0pt\it $\cal G\mskip-2mu$eometry \&\ 
$\cal T\!\!$opology $\cal P\!$ublications}}  

\def\recd{{\small Received:\qua\receiveddate\ifx\reviseddate\relax
\else\qquad Revised:\qua\reviseddate\fi\par}} 

\def\volumenumber#1{\def\thevolumenumber{#1}}
\def\volumeyear#1{\def\thevolumeyear{#1}}
\def\papernumber#1{\def\thepapernumber{#1}}
\def\pagenumbers#1#2{\def\startpage{#1}\def\finishpage{#2}}
\def\published#1{\def\publishdate{#1}}

\def\received#1{\def\receiveddate{#1}}

\def\asciiaddress#1{\def\theasciiaddress{#1}}

\long\def\asciiabstract#1{\long\def\theasciiabstract{#1}}

\let\\\par\let\thevolumenumber\relax
\let\thepapernumber\relax\let\thevolumeyear\relax\let\startpage\relax
\let\finishpage\relax\let\publishdate\relax\let\receiveddate\relax
\let\reviseddate\relax\let\theasciititle\relax
\let\theasciiauthors\relax\let\theasciiaddress\relax
\let\theasciiabstract\relax

\let\theasciiemail\relax

\ifplaintex
\font\logobig=cmssbx10 scaled 3836
\font\logomed=cmssbx10 scaled 2557
\else
\font\logobig=cmssbx10 scaled 4200
\font\logomed=cmssbx10 scaled 2800
\fi

\long\def\makeagttitle{   
\count0=\startpage
\agt\hfill      
\hbox to 45truept{\vbox to 0pt{\vglue -13truept{\logomed A\kern -.37em{\logobig 
T}\kern -.38em G}\vss}\hss}
\break
{\small Volume \thevolumenumber\ (\thevolumeyear)
\startpage--\finishpage\nl
Published: \publishdate}

\vglue .25truein

{\parskip=0pt\leftskip 0pt plus
1fil\def\\{\par\smallskip}{\Large\bf\thetitle}\par\medskip} \vglue
0.05truein

{\parskip=0pt\leftskip 0pt plus 1fil\def\\{\par}{\sc\theauthors}
\par\medskip}%
 
\vglue 0.03truein 

{\small\leftskip 25truept\rightskip 25truept{\bf Abstract}\stdspace\theabstract

{\bf AMS Classification}\stdspace\theprimaryclass
\ifx\thesecondaryclass\relax\else; \thesecondaryclass\fi\par
{\bf Keywords}\stdspace \thekeywords\par}\vglue 7truept

}   

\ifplaintex
\font\phead=cmsl9 scaled 950
\font\pnum=cmbx10 scaled 913
\font\pfoot=cmsl9 scaled 950
\headline{\vbox to 0pt{\vskip -4.5mm\line{\small\phead\ifnum
\count0=\startpage ISSN 1472-2739 (on-line) 1472-2747 (printed)
\hfill {\pnum\folio}\else\ifodd\count0\def\\{ }%
\ifx\theshorttitle\relax\thetitle\else\theshorttitle\fi\hfill{\pnum\folio}
\else\def\\{ and }{\pnum\folio}\hfill\ifx\theshortauthors\relax\theauthors
\else\theshortauthors\fi\fi\fi}\vss}}
\footline{\vbox to 0pt{\vglue 0mm\line{\small\pfoot\ifnum\count0=\startpage
\copyright\ \gtp\hfill\else
\agt, Volume \thevolumenumber\ (\thevolumeyear)\hfill\fi}\vss}}
\else
\font=cmsl9 scaled 1050
\font\lnum=cmbx10 
\font=cmsl9 scaled 1050
\makeatletter
\def\@oddhead{{\small\lhead\ifnum\count0=\startpage ISSN 1472-2739 
(on-line) 1472-2747 (printed)\hfill {\lnum\number\count0}\else\ifodd\count0
\def\\{ }\ifx\theshorttitle\relax \thetitle \else\theshorttitle\fi\hfill
{\lnum\number\count0}\else\def\\{ and }{\lnum\number\count0}
\hfill\ifx\theshortauthors\relax 
\theauthors\else\theshortauthors\fi\fi\fi}}\def\@evenhead{\@oddhead}
\def\@oddfoot{\small\lfoot\ifnum\count0=\startpage\copyright\ \gtp\hfill\else
\agt, Volume \thevolumenumber\ (\thevolumeyear)\hfill\fi}
\def\@evenfoot{\@oddfoot}
\makeatother
\fi
\let\maketitlepage\makeagttitle
\let\makeshorttitle\maketitlepage
\let\maketitle\maketitlepage

\newwrite\gtoutfile
\long\gdef\makeheadfile{  
{\def\\{, }\def\s{ }
\immediate\openout\gtoutfile head.xxx
\immediate\write\gtoutfile{To: math@arxiv.org}
\immediate\write\gtoutfile{Subject: put OR rep NNNNN:ppppp}
\immediate\write\gtoutfile{--text follows this line--}
\immediate\write\gtoutfile{Proxy-for: \ifx\theasciiauthors\relax
\theauthors\else\theasciiauthors\fi\s<\ifx\theasciiemail\relax\theemail\else\theasciiemail\fi>}
\immediate\write\gtoutfile{\noexpand\\}
\immediate\write\gtoutfile{Authors: \ifx\theasciiauthors\relax
\theauthors\else\theasciiauthors\fi}
{\def\\{ }\immediate\write\gtoutfile{Title: \ifx\theasciititle\relax
\thetitle\else\theasciititle\fi}}
\immediate\write\gtoutfile{Subj-class: GT or SG, GR etc}
\immediate\write\gtoutfile{MSC-class: \theprimaryclass\ifx\thesecondaryclass\relax\else, \thesecondaryclass\fi}
\immediate\write\gtoutfile{Journal-ref: Algebraic and Geometric Topology \thevolumenumber\s
(\thevolumeyear) \startpage-\finishpage}
\immediate\write\gtoutfile{Comments: Published in Algebraic and
Geometric Topology at}
\immediate\write\gtoutfile{\s\s\s  http://www.maths.warwick.ac.uk/agt/AGTVol\thevolumenumber/agt-\thevolumenumber-\thepapernumber.abs.html}
\immediate\write\gtoutfile{\noexpand\\}
\immediate\write\gtoutfile{}
\ifx\theasciiabstract\relax
\immediate\write\gtoutfile{\theabstract}\else
\immediate\write\gtoutfile{\theasciiabstract}\fi
\immediate\write\gtoutfile{}
\immediate\write\gtoutfile{\noexpand\\}
\immediate\write\gtoutfile{}
\immediate\closeout\gtoutfile}}  

\def\maketitlepage{\makeagttitle\makeheadfile}
\let\makeshorttitle\maketitlepage
\let\maketitle\maketitlepage

\volumenumber{1}
\volumeyear{2001}
\papernumber{1}
\published{25 October 2000}
\pagenumbers{1}{29}
\received{6 August 2000}

\usepackage{amssymb}
\usepackage{amsmath,amscd}
\usepackage{graphicx}

\newtheorem{thm}{Theorem}
\newtheorem{lem}[thm]{Lemma}

\theoremstyle{definition}

\newtheorem{defi}[thm]{Definition}

\newtheorem{remn}{Remark}

\renewcommand{\int}{\operatorname{int}}
\newcommand{\sign}{\operatorname{sign}}

\newcommand{\Z}{\mathbb{Z}}

\newcommand{\R}{\mathbb{R}}
\newcommand{\RR}{\mathcal{R}}
\newcommand{\PP}{\mathbb{P}}

\newcommand{\s}{\mathcal{S}}
\newcommand{\e}{\epsilon}

\newcommand{\imra}{\looparrowright}

\def\theenumi{\roman{enumi}}

\begin{document}

\title{Higher order intersection numbers of 2-spheres\\in 4-manifolds}
\authors{Rob Schneiderman\\Peter Teichner}
\address{Dept. of Mathematics, University of California at Berkeley\\
Berkeley, CA 94720-3840, USA}
\secondaddress{Dept. of Mathematics, University of California at
San Diego\\ La Jolla, CA 92093-0112, USA}
\email{schneido@math.berkeley.edu, teichner@euclid.ucsd.edu}

\asciiaddress{Dept. of Mathematics, University of California at Berkeley\\
Berkeley, CA 94720-3840, USA\\Dept. of Mathematics, University of California at
San Diego\\La Jolla, CA 92093-0112, USA}

\begin{abstract}

This is the beginning of an obstruction theory for deciding
whether a map $f:S^2 \to X^4$ is homotopic to a topologically
flat embedding, in the presence of fundamental group and in the
absence of dual spheres. The first obstruction is Wall's 
self-intersection number $\mu(f)$ which tells the whole
story in higher dimensions. Our second order obstruction
$\tau(f)$ is defined if $\mu(f)$ vanishes and has formally very
similar properties, except that it lies in a quotient of the
group ring of two copies of $\pi_1X$ modulo $\s_3$-symmetry
(rather then just one copy modulo $\s_2$-symmetry).
It generalizes to the non-simply connected setting the 
Kervaire-Milnor invariant defined
in \cite{FQ} and \cite{S} which corresponds to the Arf-invariant of
knots in 3-space.

We also give necessary and sufficient conditions for
moving three maps $f_1,f_2,f_3:S^2 \to X^4$ to a position in which
they have {\em disjoint} images. Again the obstruction
$\lambda(f_1,f_2,f_3)$ generalizes Wall's intersection number
$\lambda(f_1,f_2)$ which answers the same question for two
spheres but is not sufficient (in dimension~$4$) for three spheres.
In the same way as intersection numbers correspond to linking
numbers in dimension~3, our new invariant corresponds to the
Milnor invariant $\mu(1,2,3)$, generalizing the Matsumoto triple \cite{M}
to the non simply-connected setting.

\end{abstract}

\asciiabstract{ This is the beginning of an obstruction theory for
deciding whether a map f:S^2 --> X^4 is homotopic to a topologically
flat embedding, in the presence of fundamental group and in the
absence of dual spheres. The first obstruction is Wall's
self-intersection number mu(f) which tells the whole story in higher
dimensions. Our second order obstruction tau(f) is defined if mu(f)
vanishes and has formally very similar properties, except that it lies
in a quotient of the group ring of two copies of pi_1(X) modulo
S_3-symmetry (rather then just one copy modulo S_3-symmetry).  It
generalizes to the non-simply connected setting the Kervaire-Milnor
invariant which corresponds to the Arf-invariant of knots in 3-space.

We also give necessary and sufficient conditions for moving three maps
f_1,f_2,f_3:S^2 --> X^4 to a position in which they have disjoint
images. Again the obstruction lambda(f_1,f_2,f_3) generalizes Wall's
intersection number lambda(f_1,f_2) which answers the same question
for two spheres but is not sufficient (in dimension 4) for three
spheres.  In the same way as intersection numbers correspond to
linking numbers in dimension 3, our new invariant corresponds to the
Milnor invariant mu(1,2,3), generalizing the Matsumoto triple to the
non simply-connected setting.}



\maketitle

\section{Introduction} One of the keys to the success of
high-dimensional surgery theory is the following beautiful fact,
due to Whitney and Wall \cite{W}, \cite{Wh}: A smooth map
$f:S^n\to X^{2n}, n>2$, is homotopic to an embedding if and only
if a single obstruction $\mu(f)$ vanishes.  This {\em
self-intersection invariant}  takes values in a quotient of the
group ring $\Z[\pi_1X]$ by simple relations. It is defined by
observing that generically $f$ has only transverse double points
and then counting them with signs and fundamental group elements.
The relations are given by an $\s_2$-action (from changing the
order of the two sheets at a double point) and a framing
indeterminacy (from a cusp homotopy introducing a local kink).
Here $\s_k$ denotes the symmetric group on $k$ symbols.

It is well-known that the case $n=2$, $f:S^2\to X^4$, is very
different \cite{KM}, \cite{R}, \cite{HS}.
Even though $\mu(f)$ is still defined, it
only implies that the self-intersections of $f$ can be paired up
by Whitney disks. However, the Whitney moves, used in higher
dimensions to geometrically remove pairs of double points, cannot
be done out of three different reasons: The Whitney disks might
not be represented by embeddings, they might not be correctly
framed, and they might intersect $f$. Well known maneuvers on the
Whitney disks show however, that the first two conditions may
always be attained (by pushing down intersections and twisting
the boundary).

In this paper we describe the next step in an obstruction theory
for finding an embedding homotopic to $f:S^2\to X^4$ by measuring
its intersections with Whitney disks. Our main results are as
follows, assuming that $X$ is an oriented $4$-manifold.

\begin{thm}\label{thm1} If $f:S^2\to X^4$ satisfies $\mu(f) =0$
then there is a well-defined (secondary) invariant
$\tau(f)$ which depends only on the homotopy class of $f$. It takes values
in the quotient of
$\Z[\pi_1X\times\pi_1X]$ by relations additively generated by
$$
\begin{array}{crcl}
\mbox{\rm(BC)} &  (a,b) & = & -(b,a)
\\
\mbox{\rm(SC)} & (a,b) & = & -(a^{-1},ba^{-1})
\\
\mbox{\rm(FR)} & (a,1) & = & (a,a)
\\
\mbox{\rm(INT)} & (a,\lambda(f,A)) & = & (a,\omega_2(A)\cdot 1).
\end{array}
$$
where $a,b\in\pi_1X$ and $A$ represents an immersed $S^2$ or
$\R\PP^2$ in $X$. In the latter case, the group element $a$ is
the image of the nontrivial element in $\pi_1(\R\PP^2)$.
\end{thm}

If one takes the obstruction theoretic point of view seriously
then one should assume in Theorem~\ref{thm1} that in addition to
$\mu(f) =0$ all intersection numbers $\lambda(f, A) \in
\Z[\pi_1X]$ vanish. With these additional assumptions, $\tau(f)$
is defined in a quotient of $\Z[\pi_1 X \times \pi_1X]$ by
$\s_3$-and framing indeterminacies.  This is in complete analogy
with $\mu(f)$!

\begin{figure}[ht!]  
         \centerline{\includegraphics[width=0.5\hsize]{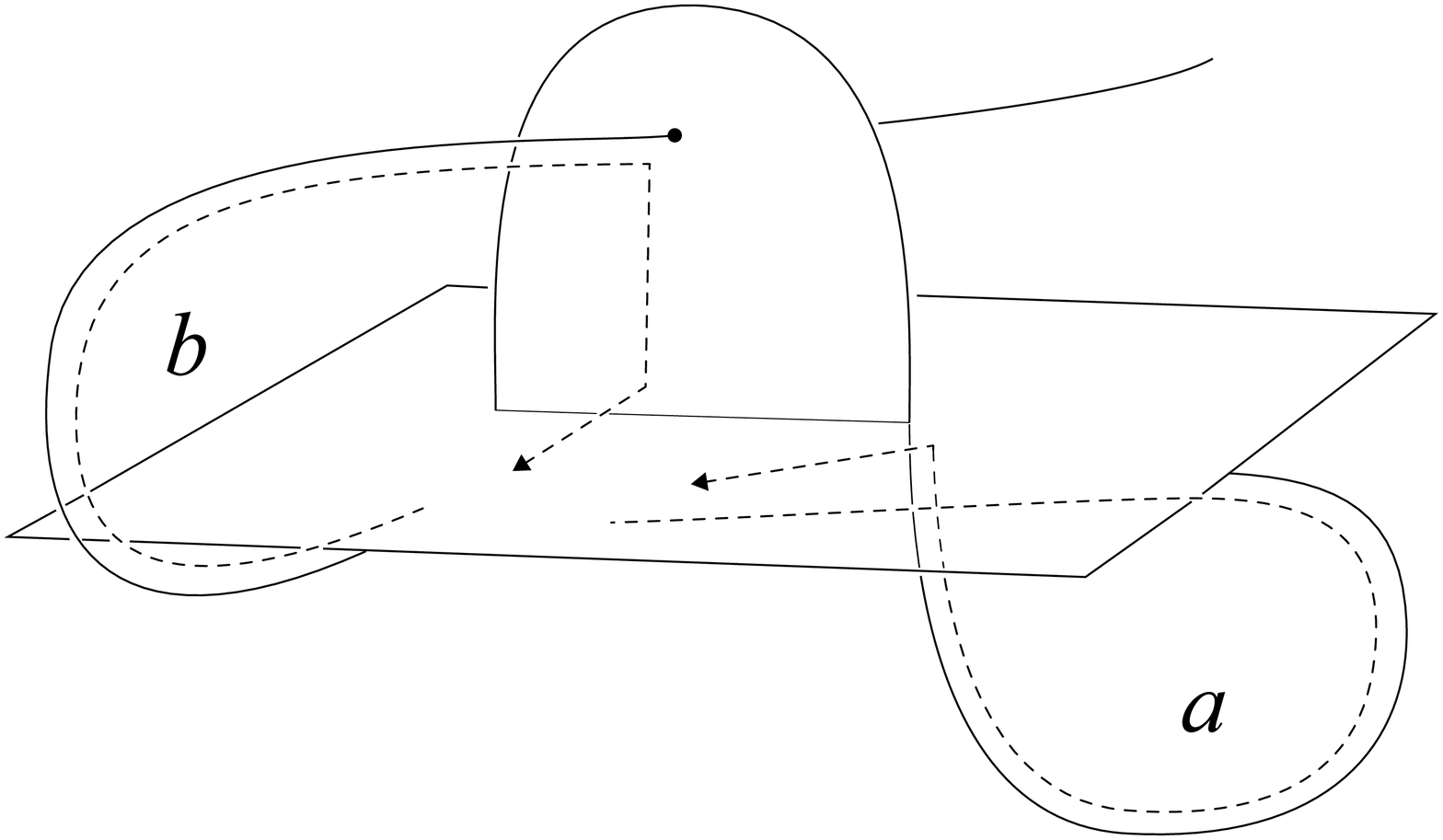}}
         \nocolon \caption{}\label{Whitneydiskab}
\end{figure}

To define $\tau(f)$, we pick framed Whitney disks for all double
points of $f$, using that $\mu(f)=0$. Then we sum up all
intersections between $f$ and the Whitney disks, recording for
each such intersection point a sign and a pair
$(a,b)\in\pi_1X\times \pi_1X$. Here $a$ measures the primary
group element of a double point of $f$ and $b$ is the secondary
group element, see Figure~\ref{Whitneydiskab}.
After introducing sign conventions, the $\s_2$-action already present
in $\mu(f)$ is easily seen to become the ``sheet change''
relation SC. The beauty of our new invariant now arises from the
fact that the other relation, which is forced on us by being able
to push around the intersections points, is of the surprisingly
simple form BC (which stands for ``boundary crossing'' of Whitney
arcs, as explained in Section~\ref{sec:tau}). In particular, this
means that the notion of {\em primary} and {\em secondary} group
elements is not at all appropriate. Moreover, one easily checks
that the two relations BC and SC together
generate an $\s_3$-symmetry on
$\Z[\pi \times \pi]$, for any group $\pi$. We will give a very
satisfying explanation of this symmetry after
Definition~\ref{def:triple}, in terms of choosing one of three
sheets that interact at a Whitney disk. The framing
indeterminacy FR comes from the being able to twist the boundary
of a Whitney disk, and the intersection relation INT must
be taken into account since one can sum a $2$-sphere into any
Whitney disk. Finally, intersections with $\R\PP^2$'s come in because of
an indeterminacy in the pairing of inverse images of double points whose
primary group elements have order two, as discovered by Stong \cite{S}.

The geometric meaning of $\tau(f)$ is
given by the following theorem, see
Section~\ref{sec:proof of thm2} for the proof.

\begin{thm}\label{thm2}  $\tau(f)=0$ if and only if $f$ is
homotopic to $f'$ such that the self-intersections of $f'$ can be
paired up by framed immersed Whitney disks with interiors
disjoint from $f'$.  In particular, $\tau(f) =0$ if $f$ is
homotopic to an embedding.
\end{thm}

The Whitney disks given by Theorem~\ref{thm2} may intersect each
other and also self-intersect. Trying to push down intersections
re-introduces intersections with $f$. Hence one expects third
(and higher) order obstructions which measure intersections among
the Whitney disks, pairing them up by secondary Whitney disks
etc.  These obstructions indeed exist in different flavors, one
has been applied in \cite{COT} to classical knot concordance.  In
a future paper we will describe obstructions living in a quotient
of the group ring of $\pi_1X \times \dots \times \pi_1X$, where
the number of factors reflects exactly the order of the
obstruction. The obstructions will be labeled by the same
uni-trivalent graphs that occur in the theory of finite type
invariants of links in $3$-manifolds. They satisfy the same
antisymmetry and Jacobi-relations as in the $3$-dimensional
setting. The reason behind this is that invariants for the {\em
uniqueness} of embeddings of $1$-manifolds into $3$-manifolds
should translate into invariants for the {\em existence} of
embeddings of $2$-manifolds into $4$-manifolds. Note that our
$\tau(f)$ corresponds the letter Y-graph and antisymmetry is
exactly our BC relation.

\begin{remn}\label{dual sphere}
It is important to note that the relation INT implies that
$\tau(f)$ vanishes in the presence of a framed dual sphere $A$.
This implies that $\tau(f)$ is not relevant in the settings of
the surgery sequence and the s-cobordism theorem. However, there
are many other settings in which dual spheres don't exist, for
example in questions concerning link concordance. The invariant
$\tau$ therefore gives a new algebraic structure on
$\pi_2(X)$ which has to be taken into account when trying to
define the correct concept of homology surgery and
$\Gamma$-groups in low dimensions.
\end{remn}
There are many examples where $\mu$ and $\lambda$ vanish but
$\tau$ is nontrivial. For infinite fundamental groups, $\tau(f)$
can take values in an infinitely generated group, see the example
in Section~\ref{sec:examples}.
If $X$ is simply-connected then $\tau(f)$ takes values in $\Z/2$
or $0$, depending on whether $f$ is spherically characteristic or
not.  In the former case, $\tau(f)$ equals the spherical
Kervaire-Milnor invariant introduced by Freedman-Quinn
\cite[Def.10.8]{FQ}. If
$X$ is closed and simply connected then $f$ has a dual sphere if
and only if it represents an indivisible class in $H_2X$. In this
case $f$ is represented by a topologically flat embedding if and
only if $\tau(f)=0$, see \cite[Thm.10.3]{FQ}. This result was
extended independently by Hambleton-Kreck \cite{HK} and
Lee-Wilczynski \cite{LW} to divisible classes $f$. They study
{\em simple} embeddings, where the fundamental group of the
complement of $f$ is abelian (and $\pi_1(X)=1$). Then there is an
additional Rohlin obstruction \cite{R} from the signature of a
certain branched cover. Moreover, these authors show that $f$ is
represented by a simple embedding if and only if $\tau(f)$ and
the Rohlin obstruction vanish. Gui-Song Li also studied the invariant
$\tau$ in a special setting \cite{L} which actually motivated our
discussion.

The $\s_3$-symmetry of $\tau$ comes from the fact that we cannot
distinguish the three sheets interacting at a Whitney disk. It is
therefore not surprising that there is a simpler version of this
invariant, defined for three maps $f_1,f_2,f_3:S^2 \to X^4$. It
can be best formulated by first identifying $\pi\times\pi$ with
the quotient $(\pi\times\pi \times \pi)/\Delta(\pi)$, where
$\Delta$ denotes the diagonal {\em right} action of $\pi$. Let $
\Lambda:= \Z[(\pi\times\pi \times \pi)/\Delta(\pi)]$ which is a
$\Z[\pi\times\pi \times \pi]$-module via left multiplication. It
also has an obvious $\s_3$-action by permuting the three factors.
This action agrees with the action generated by BC and SC if we
make the correct identification of $ \Lambda$ with
$\Z[\pi\times\pi]$.

   Now define
$\RR$ to be the $\Z[\pi\times\pi \times \pi]$-submodule of $
\Lambda$ generated by $$ (a,b, \lambda(f_3,A)), (a,
\lambda(f_2,A),b), (\lambda(f_1,A),a,b)\in \Lambda $$ where
$a,b\in\pi:=\pi_1X$ and $A\in\pi_2X$ are arbitrary. The following
result will be proven in Section~\ref{sec:triples}.

\begin{thm}\label{triples}
In the above notation, assume that $ \lambda(f_i,f_j)=0$ for
$i\neq j$. Then there is a well defined secondary obstruction $$
\lambda(f_1,f_2,f_3)\in \Lambda/\RR $$ which only depends on the
homotopy classes of the $f_i$. It vanishes if and only if the
$f_i$ are homotopic to three maps with disjoint images. Moreover,
it satisfies the following algebraic properties (where
$a,a'\in\pi$ and $\sigma\in\s_3$):
\begin{enumerate}
\item $\lambda(a\cdot f_1 + a'\cdot f_1', f_2, f_3)= (a,1,1)\cdot
\lambda(f_1,f_2,f_3) + (a',1,1)\cdot \lambda(f_1',f_2,f_3)$,
\item $ \lambda(f_{ \sigma(1)}, f_{\sigma(2)}, f_{ \sigma(3)}) =
\lambda(f_1,f_2,f_3)^\sigma$,
\item  $\lambda(f,f,f) = \sum_{\sigma\in \s_3} \tau(f)^\sigma$ if $f$ has
trivial normal bundle,
\item $\tau(f_1+f_2+f_3)-\tau(f_1+f_2)-\tau(f_1+f_3)-\tau(f_2+f_3)
+\tau(f_1)+\tau(f_2)+\tau(f_3) = \lambda(f_1,f_2,f_3)$,
\item $\tau(a\cdot f)=a^{-1} \tau(f)a$.
\end{enumerate}
\end{thm}
These properties are the precise analogues of the fact that
Wall's $(\lambda,\mu)$ is a ``quadratic form'' on $\pi_2(X)$ (or
a hermitian form with a quadratic refinement), see \cite[\S 5]{W}.
To make this precise one has to identify $\Z[(\pi \times
\pi)/\Delta(\pi)]$ with $\Z[\pi]$ via the map
$$
(a,b)\mapsto a\cdot b^{-1}.
$$
Then the usual involution $a\mapsto a^{-1}$ corresponds to flipping the
two factors which explains why (ii) above generalizes the notion of a
{\em hermitian} form.
It would be nice if one could formalize these ``cubic forms'', guided by
the above properties.

Note that if the primary intersection numbers $\lambda(f_i,A)$
vanish for all $A\in\pi_2X$, then $\lambda(f_1, f_2, f_3)$ is
well defined as an element of $\Z[\pi_1X\times\pi_1X]$.
     If $X$ is simply-connected, this reduces to
the Matsumoto triple from \cite{M}.
Garoufalidis and Levine have also introduced equivariant
$\mu(1,2,3)$-invariants in \cite{GL} for null homotopic circles in a
$3$-manifold $N^3$. These invariants agree with our triple $\lambda$
applied to three disks in $N^3 \times [0,1]$ that display the null
homotopies. If one wants to get spheres instead of disks, one should
attach 2-handles to all the circles, and glue the cores to the null
homotopies. In
\cite{GL} the indeterminacies of the invariants are not discussed but it
is shown that they agree for two links if and only if they
are surgery equivalent \cite[Thm.5]{GL}.

The obstruction $\lambda(f_1,f_2,f_3)$ generalizes Wall's
intersection number $\lambda(f_1,f_2)$ which answers the disjointness
question for two spheres
    but is not sufficient (in dimension~$4$) for three spheres. In an
upcoming paper we will describe necessary and sufficient
obstructions for making $n$ maps $f_1,\dots,f_n:S^2\to X^4$
disjoint. The last obstructions will lie in the group ring of
$(n-1)$ copies of $\pi_1X$, assuming that all previous
obstructions vanish.

The current paper finishes with Section~\ref{sec:general tau} by
giving the following generalization of Theorem~\ref{thm1} and
Theorem~\ref{thm2} to the case of arbitrarily many maps.
\begin{thm}\label{n spheres}  Given $f_1, \dots, f_n: S^2\to X^4$ with
vanishing primary $\mu$ and $\lambda$- obstructions, there exists
a well-defined secondary obstruction $\tau(f_1,\ldots, f_n)$
which depends only on the homotopy classes of the $f_i$.  This
invariant  vanishes if and only if, after a homotopy, all
intersections and self-intersections can be paired up by Whitney disks
with interiors disjoint from all $f_i(S^2)$.  This is in
particular the case if the $f_i$ are homotopic to disjoint
embeddings.
\end{thm}
The invariant $\tau(f_1,\ldots,f_n)$ takes values in a quotient
of ${n \choose 1}+2{n \choose 2}+{n \choose 3}$ copies of
$\Lambda$ which reflects the number of different combinations of
possible intersections between Whitney disks and spheres.

In this paper we assume that our 4-manifolds are oriented and we
work in the smooth category.  However, our methods do not
distinguish the smooth from the topological category since the
basic results on topological immersions \cite{FQ} imply a
generalization of our results to the topological world.

\section{Preliminaries}\label{sec:prelim}

We refer the reader to the book by Freedman and Quinn \cite[\S
1]{FQ} for the basic definitions of things like Whitney disks,
Whitney moves, finger moves and Wall's intersection and
self-intersection numbers $ \lambda $ and $ \mu$ (see also
\cite{W}). We only make a couple of summarizing remarks.

Let $f:S^2\to X^4$ be a smooth map. After a small perturbation we
may assume that $f$ is a generic immersion. This means that the
singularities of $f$ consist only of transverse self-intersection
points. Furthermore, we may perform some cusp homotopies to get
the signed sum of the self-intersection points of $f$ to be zero
as an integer.
    By an old theorem of Whitney, immersions $f:S^2 \imra X^4$ as
above, modulo regular homotopy, are the same as homotopy classes
of maps $S^2\to X^4$. We will thus assume that our maps $S^2\imra
X^4$ are immersions with only transverse self-intersections whose
signed sum is zero. Then we work modulo regular homotopy. The
advantage of this approach comes from the fact that by general
position, a regular homotopy is (up to isotopy) the composition
of finitely many finger moves and then finitely many Whitney
moves. This implies that $\mu(f)$ is well-defined in the quotient
of the group ring $\Z[\pi_1X]$ by the $\s_2$-action $a\mapsto
a^{-1}$.

Let $f:S^2\imra X^4$ be a generic immersion and let $p,q\in X$ be
double points of opposite sign. Choose two embedded arcs in $S^2$
connecting the inverse images of $p$ to the inverse images of $q$
but missing each other and all other double points of $f$. The
image $\gamma$ of the union of these {\em Whitney arcs} is called
a {\em Whitney circle} for $p,q$ in $X$. Let $W:D^2 \imra X$ be
an immersion which is an embedding on the boundary with
$W(S^1)=\gamma$. The normal bundle of $W $ restricted to $\gamma$
has a canonical nonvanishing section $s_\gamma$ which is given by
pushing $\gamma$ tangentially to $f$ along one of the Whitney
arcs and normally along the other. Therefore, the relative Euler
number of the normal bundle of $W $ is a welldefined integer. If
one changes $W$ by a (nonregular) cusp homotopy then the Euler
number changes by $\pm 2$,
    see \cite[\S 1.3]{FQ}. This implies that one really has a
$\Z/2$-valued framing invariant. An additional boundary twist can
be used to change the Euler number by one. Note that this
introduces an intersection between $W$ and $f$ and thus does not
preserve the last property below.
\begin{defi} \label{def:Whitney-disk}
Let $W:D^2 \imra X$ be an immersion as above.
\begin{enumerate}
\item If $W$ has vanishing relative Euler number, then it is
called a {\em framed Whitney disk}. Some authors also add the
adjective {\em immersed} but we suppress it from our notation.
\item If in addition $W$ is an embedding with interior
disjoint from $f$, then $W$ is called an {\it embedded Whitney disk} for $f$.
\end{enumerate}
\end{defi}
If $W$ is an embedded Whitney disk one can do the Whitney move to
remove the two double points $p$ and $q$. If one of the
conditions for an embedded Whitney disk are not satisfied the
Whitney move can still be done but it introduces new
self-intersections of $f$.

The vanishing of $\mu(f)$ means that the
   double points of $f$ occur in {\it cancelling pairs}
   with opposite signs and contractible Whitney circles. Therefore,
there exists a collection of framed Whitney disks pairing up all
the double points of $f$.

\section{The invariant $\tau$} \label{sec:tau}

    In this section we define the invariant
$\tau$ of Theorem~\ref{thm1} in terms of fundamental group
elements which are determined by two kinds of intersections:
Intersections between the interiors of Whitney disks and the
sphere $f$ and intersections among the boundary arcs of the
Whitney disks.
    The definition will involve first making choices and then
modding out the resulting indeterminacies. Many of these
indeterminacies will be noted during the course of the definition
but a complete proof that $\tau(f)$ is indeed well-defined (and
only depends on the homotopy class of $f$) will be given in
Section~\ref{sec:tau well-defined}.

In the following discussion we will not make a distinction
between $f$ and its image unless necessary. Also, basepoints and
their connecting arcs ({\it whiskers}) will be suppressed from
notation.

   Let $f: S^2 \imra X^4$ be an oriented generic
immersion with vanishing Wall self intersection invariant
   $\mu(f)=0$. As explained above, the vanishing of $\mu(f)$ implies that we may
choose framed
   Whitney disks $W_i$ for all canceling pairs
$(p^{+}_{i},p^{-}_{i})$ of double points of $f$ where $\sign
(p^{+}_{i}) =+1=-\sign (p^{-}_{i})$. We may assume that the
interiors of the Whitney disks are transverse to $f$. The
boundary arcs of the $W_i$ are allowed to have transverse
intersections and self-intersections (as arcs in $f$).

\begin{remn}\label{weak disks}
In the literature it is often assumed that a collection of
Whitney disks will have disjointly embedded boundary arcs.
Whitney disks with immersed boundaries were called ``weak''
Whitney disks in \cite{FQ} and \cite{S}. Allowing such weak
Whitney disks in the present setting will simplify the proof that
$\tau$ is well-defined.
\end{remn}

\begin{figure}[ht!] 
         \centerline{\includegraphics[width=0.6\hsize]{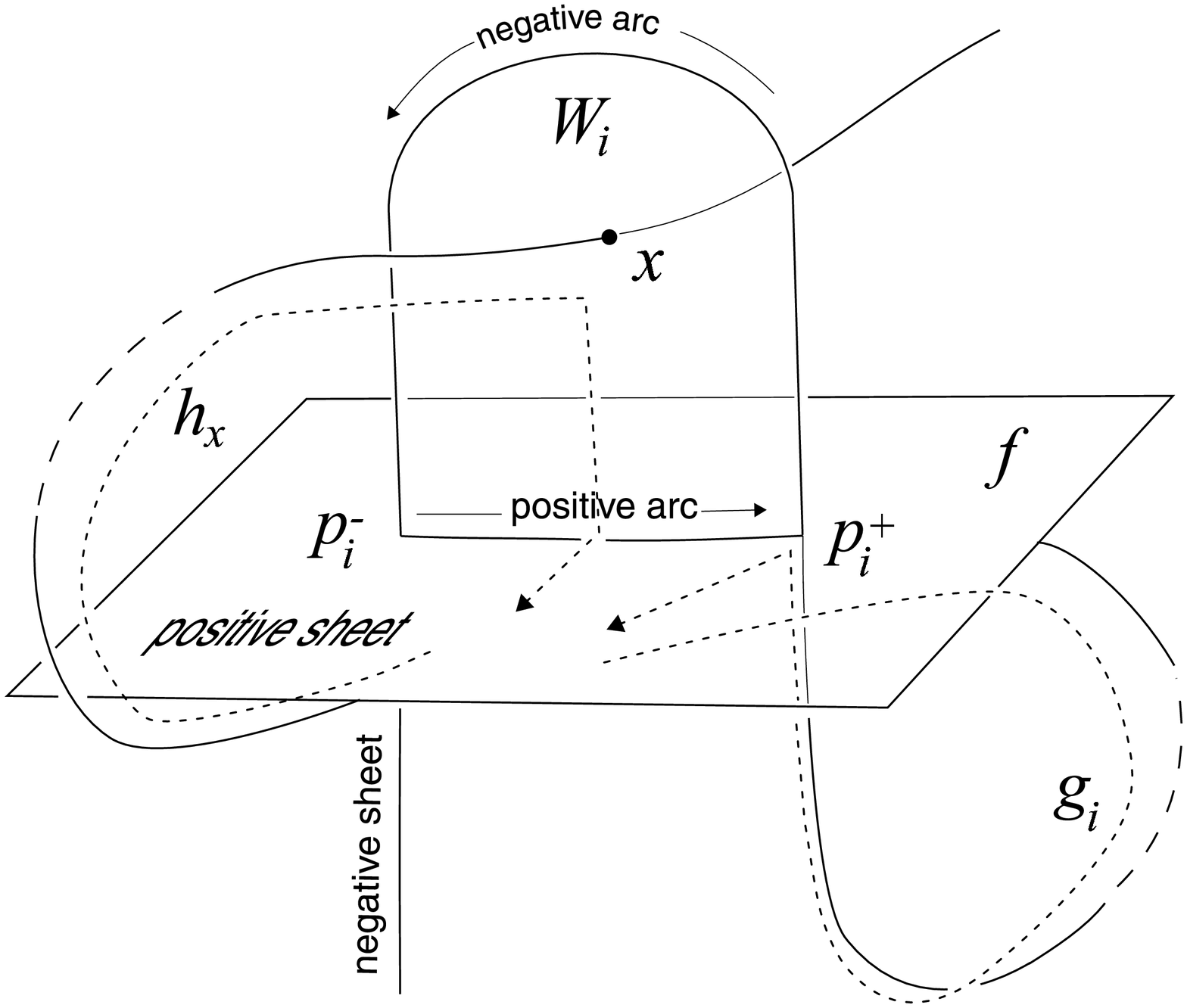}}
         \caption{Whitney disk conventions.}
         \label{Whitneydisk conventions}

\end{figure}

   For each $W_i$ choose a preferred arc of $\partial W_i$
which runs between $p^{+}_{i}$ and $p^{-}_{i}$. We will call this
chosen arc the \emph{positive arc} of $W_i$ and the arc of
$\partial W_i$ lying in the other sheet will be called the
\emph{negative arc}. We will also refer to a neighborhood in $f$
of the positive (resp.\ negative) arc as the \emph{positive}
(resp.\ \emph{negative}) \emph{sheet} of $f$ near $W_i$. This
choice of positive arc determines an orientation of $W_i$ as
follows: Orient $\partial W_i$ from $p^{-}_{i}$ to $p^{+}_{i}$
along the positive arc and back to $p^{-}_{i}$ along the negative
arc. The positive tangent to $\partial W_i$ together with an
outward pointing second vector orient $W_i$. The choice of
positive arc also determines a fundamental group element $g_i$ by
orienting the double point loops to change from the negative
sheet to the positive sheet at the double points $p^{\pm}_i$.
(See Figure~\ref{Whitneydisk conventions})
Note that changing the choice of positive arc
reverses the orientation of $W_i$ and changes $g_i$ to
$g_i^{-1}$. These orientation conventions will be assumed in the
definitions that follow.

For a point $x\in \int W_i \cap f$, define $h_x \in \pi_1 (X)$
from the following loop: Go first along $f$ from the basepoint to
$x$, then along $W_i$ to the positive arc of $W_i$, then back to
the basepoint along $f$. This loop (together with the whisker on
$f$) determines $h_x$. (See Figure~\ref{Whitneydisk
conventions}). Note that changing the
choice of positive arc for $W_i$ changes $h_x$ to
$h_x g_{i}^{-1}$.

\textsc{Notation convention}: For a sum of elements in
the integral group ring $\Z [\pi_1 X \times \pi_1 X]$
with a common first component it will
sometimes be convenient to write the sum inside the parentheses:
$$ (g, {\textstyle \sum}_j n_j\cdot g_j ) := {\textstyle \sum}_j
n_j\cdot(g,g_j ) \in \Z [\pi_1 X \times \pi_1 X]. $$

We can now begin to measure
intersections between the Whitney disks and $f$ by defining
$$
I(W_i):=(g_i ,{\textstyle \sum}_{x} \sign (x) h_x) \in \Z [\pi_1 X
\times \pi_1 X]
$$
where the sum is over all  $x\in \int W_i \cap f$, and $\sign
(x)=\pm 1$ comes from the orientations of $f$ and $W_i$ as above.

\begin{remn} \label{rem:Wall-type}$I(W_i)$ encodes Wall-type
intersections between
$W_i$ and $f$ and can be roughly written as
$(g_i,\lambda(W_i,f))$. We will see that $\tau(f)$ measures to
what extent the \emph{sum} over $i$ is well-defined. This idea
will be developed further in Section~\ref{sec:triples}.
\end{remn}

Next we set up notation to measure intersections between the
boundaries of the Whitney disks. Denote the positive arc (resp.
negative arc) of $W_i$ by $\partial_{+}{W_i}$ (resp.
$\partial_{-}{W_i}$).
Let $y$ be any point in
$\partial_{\e_{i}}
W_{i}\cap\partial_{\e_{j}}W_{j}$
where the ordered basis
$(\overrightarrow{{\partial W_{i}}},\overrightarrow{{\partial
W_{j}}})$
agrees with the orientation of $f$ at $y$.
Define
$$
J(y):=\e_{i}\e_{j}(g_{i}^{\e_{i}},g_{j}^{\e_{j}}) \in \Z [\pi_1 X
\times \pi_1 X]
$$
where
$\e_k \in \{+,- \}\widehat{=} \{+1,-1 \}$.

\begin{figure}[ht!] 
       \cl{\includegraphics[width=.8\hsize]{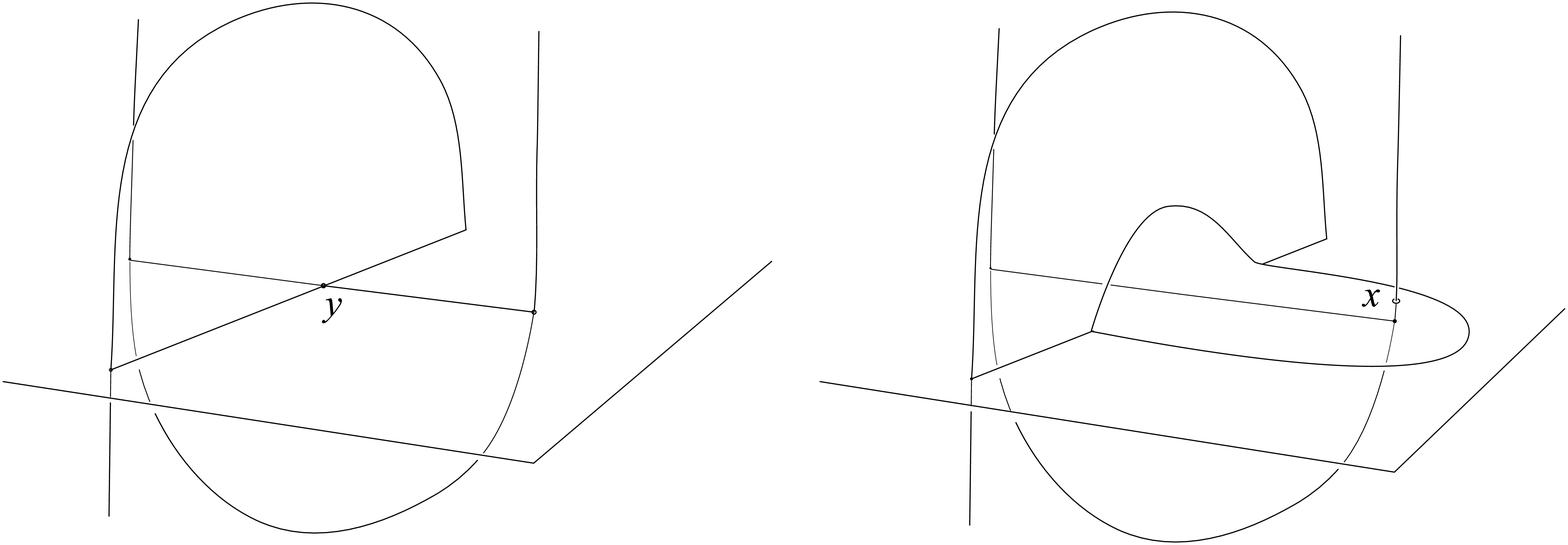}}
         \caption{Eliminating an intersection between
         Whitneydisk boundaries creates an interior intersection
         between a Whitneydisk and $f$.}
         \label{y-push}

\end{figure}

Note that by pushing $W_i$ along $\partial W_j$, as in
Figure~\ref{y-push}, $y$ could be eliminated at the cost of
creating a new intersection point $x \in \int W_i \cap f$ with
$h_x=g_{j}^{\e_{j}}$ whose contribution to $I(W_i)$ would be
$\e_{i}\e_{j}(g_{i}^{\e_{i}},g_{j}^{\e_{j}})=J(y)$. Similarly,
$y$ could also be eliminated by pushing $W_j$ along $\partial W_i$
which would create a new intersection point in $\int W_j
\cap f$; however this new intersection point would contribute
$-\e_{i}\e_{j}(g_{j}^{\e_{j}},g_{i}^{\e_{i}})$ to $I(W_j)$
illustrating the need for the BC relation.

Having made the above choices we now define our invariant:
\begin{defi}\label{tau-defi}
For $f$ as above, define
$$ \tau(f):=\sum_{i}I(W_i)+\sum_y J(y)\in {\Z [\pi_1 X \times \pi_1
X]}/{\RR} $$
where the first sum is over all Whitney disks and the second sum
is over all intersections between the boundaries of the Whitney disks.
\end{defi}

The relations $\RR$ are additively generated by the following
equations:
$$
\begin{array}{crcl}
\mbox{(BC)} \ \ &  (a,b) & = & -(b,a)
\\
\mbox{(SC)} \ \ & (a,b) & = & -(a^{-1},ba^{-1})
\\
\mbox{(FR)} \ \ & (a,1) & = & (a,a)
\\
\mbox{(INT)} \ \ &(a,\lambda(f,A))& = & (a,\omega_2(A)\cdot 1).
\end{array}
$$

Here $a$ and $b$ are any elements in $\pi_1X$ and $1\in \pi_1X$
is the trivial element. The labels BC, SC, FR and INT stand for
``boundary crossing'', ``sheet change'', ``framing'' and
``intersections'', respectively. As discussed above, the BC
relation comes from the indeterminacy in the $J$-component of
$\tau(f)$ and the other three relations come from indeterminacies
in the $I$-component of $\tau(f)$. The sheet change SC has been
already discussed above, whereas FR comes about as follows:
Changing a Whitney disk $W_i$ by a boundary twist around the
positive (resp.\ negative) arc creates $x\in \int W_i \cap f$ with
$h_x=1\in \pi_1 (X)$ (resp.\ $h_x=g_i \in \pi_1X$). After
introducing an even number of boundary twists, the correct
framing on $W_i$ can be recovered by introducing interior twists
(if necessary); this changes $I(W_i)$ by $(g_i,n +m\cdot g_i)$
where $n$ and $m$ are integers and $n \equiv m$ modulo 2. Note
that by the BC and SC relations we have $$
(a,1)=-(1,a)=(1,a)=-(a,1) \quad\Longleftrightarrow \quad (a,2)=0
$$ and hence the relation FR above is all that is needed in
addition to this relation.

The INT relation comes from changing the homotopy class of $W_i$
by tubing into any 2-sphere $A$. After correcting the framing on
$W_i$ by boundary-twists (if necessary) this changes $I(W_i)$ by
$(g_i,\lambda(f,A)+\omega_2(A)\cdot 1)$. The $\omega_2$ term is
only defined modulo 2 but still makes sense in ${\Z [\pi_1 X \times
\pi_1 X]}/{\RR}$
because $(g_i,2)=0$.

    The INT relation should, in fact, be interpreted in a more
general sense which we now describe. This goes back to an error
in \cite{FQ} as corrected by Stong \cite{S}. In the case where
$a\in\pi_1X$ satisfies $a^{2}=1$, then we allow $A$ to be not
just any immersed 2-sphere in $X$ but also any immersed $\R\PP^2$
in $X$ representing $a$, that is, $a$ is the image of the
generator of the fundamental group of $\R\PP^2$. In general, a
Wall intersection between an immersed $\R\PP^2$ and $f$ is not
well-defined because $\R\PP^2$ is not simply connected. However,
the expression $(a,\lambda(f,A))$ makes sense in
${\Z [\pi_1 X \times \pi_1 X]}/{\RR}$ because of the SC relation which accounts
exactly for the fundamental group of $\R\PP^2$ and the
orientation-reversing property of any non-trivial loop. As will
be seen in the proof of Theorem~\ref{thm1} below, the INT
relation in this case corresponds to a subtle indeterminacy in
the choice of Whitney disk for a cancelling pair of double points
whose group element $a$ has order~two.

\begin{remn}\label{rem:J=0}
It is interesting to note that one can always use finger moves to
eliminate all intersections between $f$ and the interiors of the
$W_i$ so that $\tau(f)$ is given completely in terms of the $J$
contributions from intersections between the boundary arcs
$\partial{W_i}$. On the other hand, the boundary arcs
$\partial{W_i}$ can always be made to be disjointly embedded
(Figure~\ref{y-push}) so that $\tau(f)$ is completely given in terms of the
contributions to the $I(W_i)$ coming from intersections between
$f$ and the interiors of the $W_i$.
\end{remn}

\begin{remn}\label{rem:1-conn}
If $X$ is simply connected then ${\Z [\pi_1 X \times \pi_1 X]}/{\RR}$
is $\Z/2$ or
$0$ depending on whether $f$ is spherically characteristic or
not. Moreover, $\tau(f)$ reduces to the spherical Kervaire-Milnor
invariant $km(f)\in \Z_2$  described in \cite{FQ} and \cite{S}.
If $X$ is not simply connected then $km(f)$ is equal to $\tau(f)$
mapped forward via $\pi_1X \to \{1\}$.
\end{remn}

\begin{remn}\label{unframed} One can modify the relations $\RR$ to get a
version of $\tau$ that ignores the framings on the Whitney disks
and a corresponding unframed version of Theorem~\ref{thm2}: Just
change the FR relation to $(a,1)=0$ and note that this kills the
$\omega_2$ term in the INT relation.
\end{remn}

\section{Examples}\label{sec:examples}

In this section we describe examples of immersed spheres
$f:S^2\looparrowright X$ such that ${\Z [\pi_1 X \times \pi_1 X]}/{\RR}$ is
infinitely generated and $\tau(f)$ realizes any value in
${\Z [\pi_1 X \times \pi_1 X]}/{\RR}$.

\begin{figure}[ht!] 
         \centerline{\includegraphics[width=0.7\hsize]{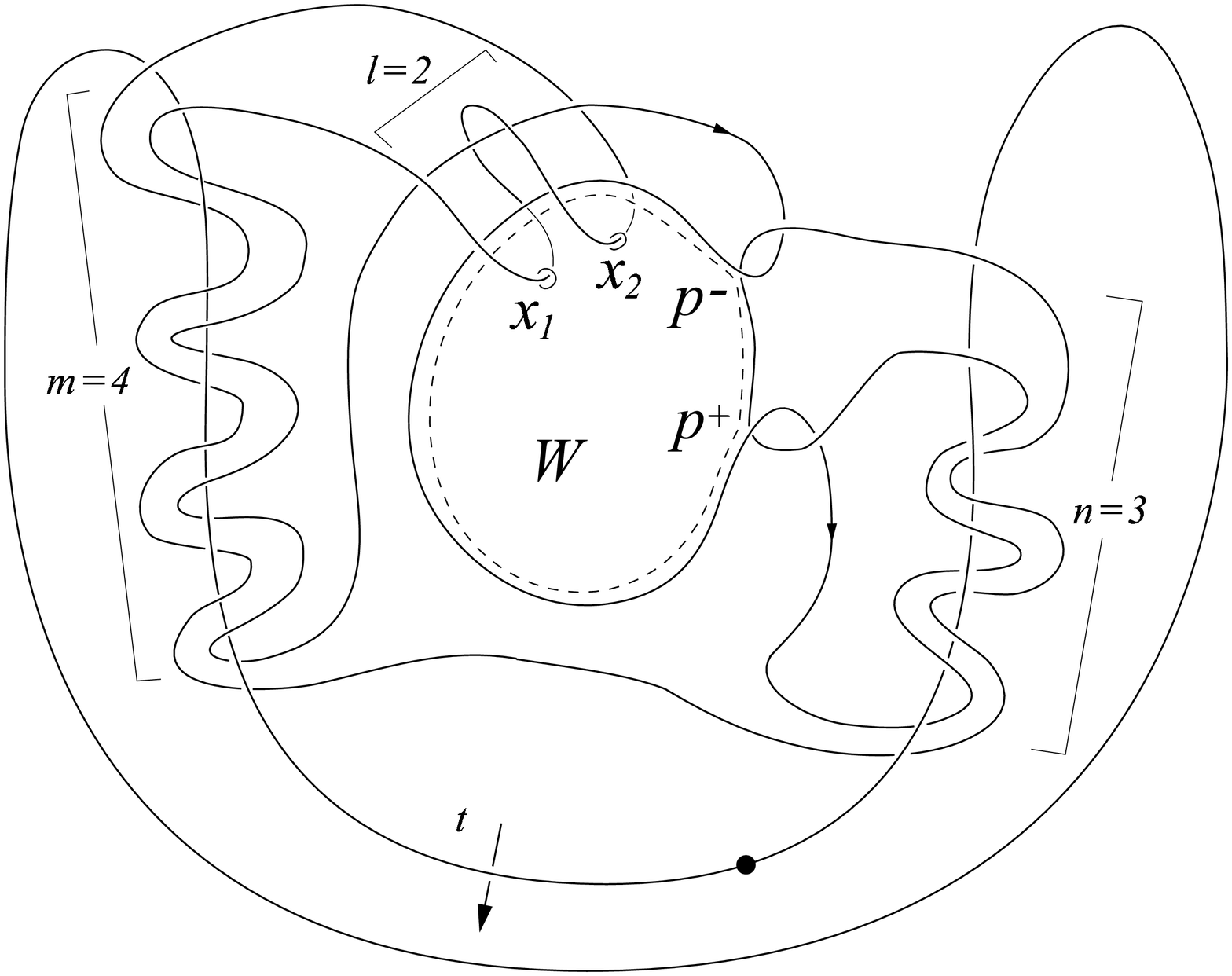}}
         \nocolon \caption{}\label{example}

\end{figure}

   Figure~\ref{example} shows the case $(l,m,n)=(2,4,3)$ of a
family of 2-component links in $S^3=\partial B^4$ indexed by
triples of integers. A 4-manifold $X$ is described by removing a
tubular neighborhood of the obvious spanning disk (pushed into
$B^4$) for the dotted component and attaching a 0-framed 2-handle
to the other component. A meridian $t$ to the dotted component
generates $\pi_1X=\langle t \rangle \cong \Z$. The other
component is an ``equator'' to an immersed 2-sphere
$f:S^2\looparrowright X$ with $\mu(f)=0$ which generates
$\pi_2(X)$ as we now describe. One hemisphere of $f$ is the core
of the 2-handle. The other hemisphere of $f$ is the trace of a
nullhomotopy of the equator in a collar of $X$. This nullhomotopy
is described by changing the two crossings labeled $p^\pm$ and
then capping off the resulting unknot with an embedded disk. The
only two double points of $f$ come from the crossing changes of
the nullhomotopy and form a canceling pair with corresponding
group element $t^n$. The dashed loop indicates a collar of a
framed embedded Whitney disk $W$ for this cancelling pair. The
interior of $W$ intersects $f$ in $l$ points $x_j$,
$j=1,2,\ldots,l$, with $h_{x_j}=t^m$ for all $j$. It follows that
$\tau(f)=l(t^n,t^m)\in\Z [\pi_1 X \times \pi_1 X]$. By band summing
different members
of this family of links one can generalize this construction to
describe $f:S^2\looparrowright X$ with $\mu(f)=0$, $\pi_1X=\Z$,
and $\pi_2(X)=\langle f \rangle$ such that $\tau(f)$ realizes any
value in ${\Z [\pi_1 X \times \pi_1 X]}/{\RR}$.

Since $f$ generates $\pi_2(X)$ and $\mu(f)=0$ (and
$\omega_2(f)=0$) the INT relation is trivial. So in this case
$\Z [\pi_1 X \times \pi_1 X]/{\RR}$ is the quotient of $\Z[\Z^2]$ by
the order 6
orbits of the $\s_3$ action generated by the SC and BC relations
together with the identification of the diagonal with the first
factor given by the FR relation. In particular ${\Z [\pi_1 X \times
\pi_1 X]}/{\RR}$
is not finitely generated.

\section{Proof of Theorem~\ref{thm1}}\label{sec:tau well-defined}

To prove Theorem~\ref{thm1} we first show that $\tau(f)$ (as
defined in Section~\ref{sec:tau}) is well-defined by considering
all the possible indeterminacies in the Whitney disk construction
used to define $\tau$ and then check that $\tau(f)$ is unchanged
by finger moves and Whitney moves on $f$ which generate
homotopies of $f$. The outline of our proof mirrors the arguments
in \cite{FQ}, \cite{S} with the added complications of working
with signs and $\pi_1X$.

In the setting of Section~\ref{sec:tau}, let $f:S^{2}\imra X$ be
a generic immersion with $\mu(f)=0$ and cancelling pairs of
double points $(p^{+}_{i},p^{-}_{i})$ paired by framed Whitney disks $W_i$
with chosen positive arcs.

Changing the choice of positive arc for a Whitney disk $W_i$
changes the orientation of $W_i$ and changes the contribution to
$I(W_i)$ of each $x\in \int(W_i)\cap f$ from $\pm (g_i,h_x)$ to
$\mp (g_{i}^{-1},h_x g_{i}^{-1})$. This does not change $\tau(f)$
by the SC relation.

Consider the effect on $I(W_i)$ of changing the interior of a
Whitney disk $W_i$: Let  ${W_i}'$ be another framed Whitney disk
with $\partial {W_i}'=\partial {W_i}$. After performing boundary
twists on $W_i$ (if necessary), $W_i$ (minus a small collar on
the boundary) and ${W_i}'$ (with the opposite orientation and
minus a small collar on the boundary) can be glued together to
form an immersed 2-sphere $A$ which is transverse to $f$. If $n$
boundary twists were done around the positive arc and $m$
boundary twists were done around the negative arc we have $$
I(W_i)-I({W_i}')=(g_i,\lambda(f,A))+ n(g_i,1) + m(g_i,g_i)=
(g_i, \lambda(f,A)+(n+m)\cdot 1) $$ where the second equality
comes from the FR relation. Since (before the boundary twists)
$W_i$ and ${W_i}'$ were correctly framed we have $n+m \equiv
\omega_2(A)$ mod 2. It follows that $ I(W_i)-I({W_i}')$ equals
zero in ${\Z [\pi_1 X \times \pi_1 X]}/{\RR}$ by the INT relation. Thus the
contribution of $I(W_i)$ to $\tau(f)$ only depends on $\partial
W_i$.

\begin{figure}[ht!] 
         \centerline{\includegraphics[width=0.5\hsize]{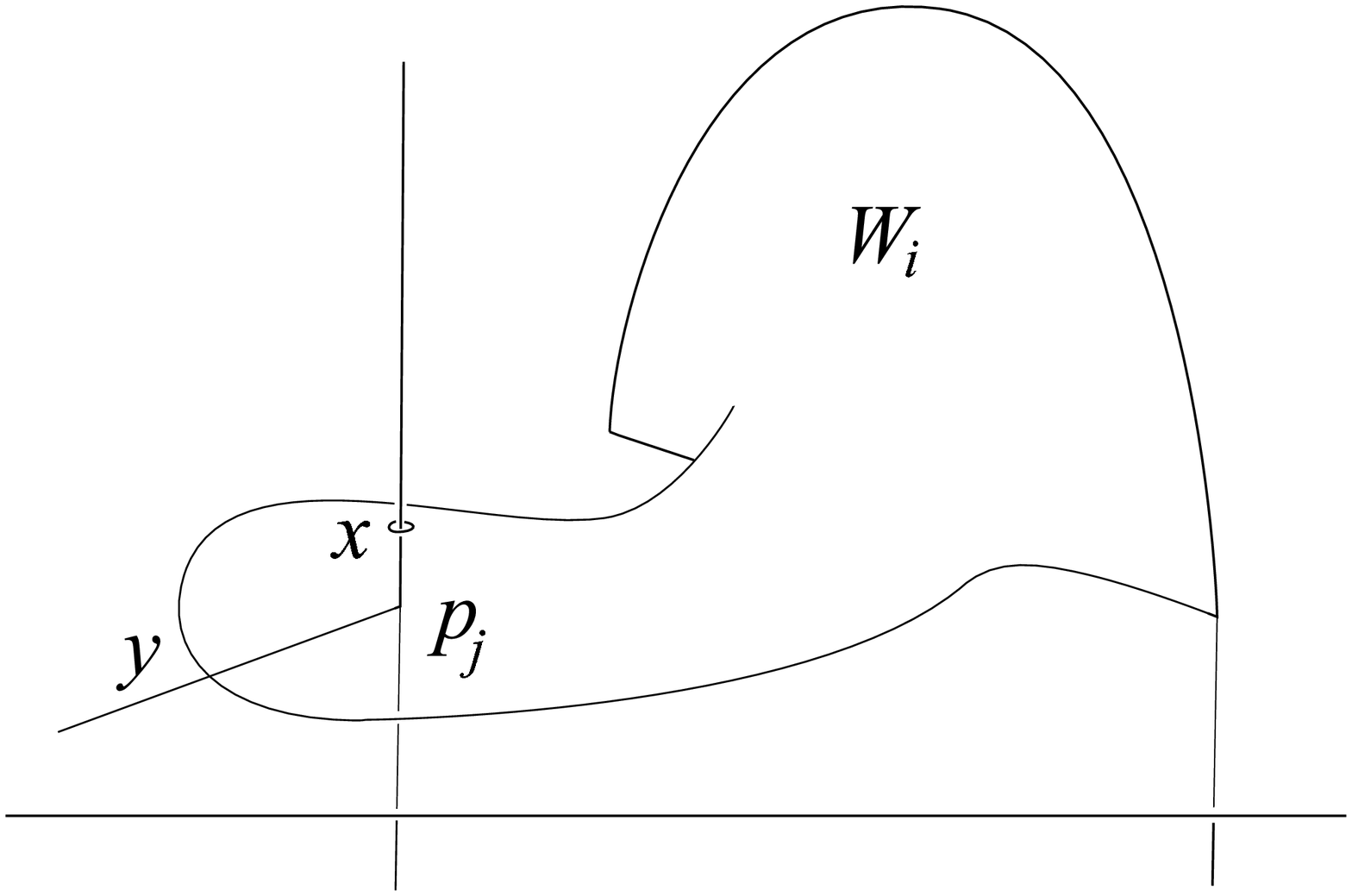}}
         \nocolon \caption{}\label{pushacross}

\end{figure}

Now consider changing $\partial W_i$ by a regular homotopy rel
$(p^{+}_{i},p^{-}_{i})$. Such a homotopy extends to a regular
homotopy of $W_i$ which is supported in a small collar on
$\partial W_i$. Away from the double points of $f$ the homotopy
can create or eliminate pairs of intersections between boundary
arcs. These pairs have canceling $J$ contributions so that
$\tau(f)$ is unchanged. When the homotopy crosses a double point
$p_j$ of $f$ a new intersection $x\in f\cap \int(W_i)$ and a new
intersection $y\in \partial_{\e_{i}}
W_{i}\cap\partial_{\e_{j}}W_{j}$ are created
(see Figure~\ref{pushacross}). One can check that
the contribution of $x$ to $I(W_i)$ is cancelled in
${\Z [\pi_1 X \times \pi_1 X]}/{\RR}$
by $J(y)$: If $\e_{i}=+=\e_{j}$ and the orientation
of $f$ at $y$ agrees with
$(\overrightarrow{{\partial W_{i}}},\overrightarrow{{\partial W_{j}}})$
then $x$ contributes
$-(g_i,g_j)$ and $J(y)=(g_i,g_j)$; If
$\e_{i}=+=\e_{j}$ and the orientation
of $f$ at $y$ agrees with
$(\overrightarrow{{\partial W_{j}}},\overrightarrow{{\partial W_{i}}})$
then $x$ contributes
$+(g_i,g_j)$ and $J(y)=(g_j,g_i)=-(g_i,g_j)$ by the BC relation.
Other cases are checked similarly.
Since any two collections of
immersed arcs (with the same endpoints) in a 2-sphere are
regularly homotopic (rel $\partial$), it follows that $\tau(f)$
does not depend on the choices of Whitney disks for given
pairings $(p^{+}_{i},p^{-}_{i})$ of the double points of $f$.

\begin{figure}[ht!] 
         \centerline{\includegraphics[width=0.5\hsize]{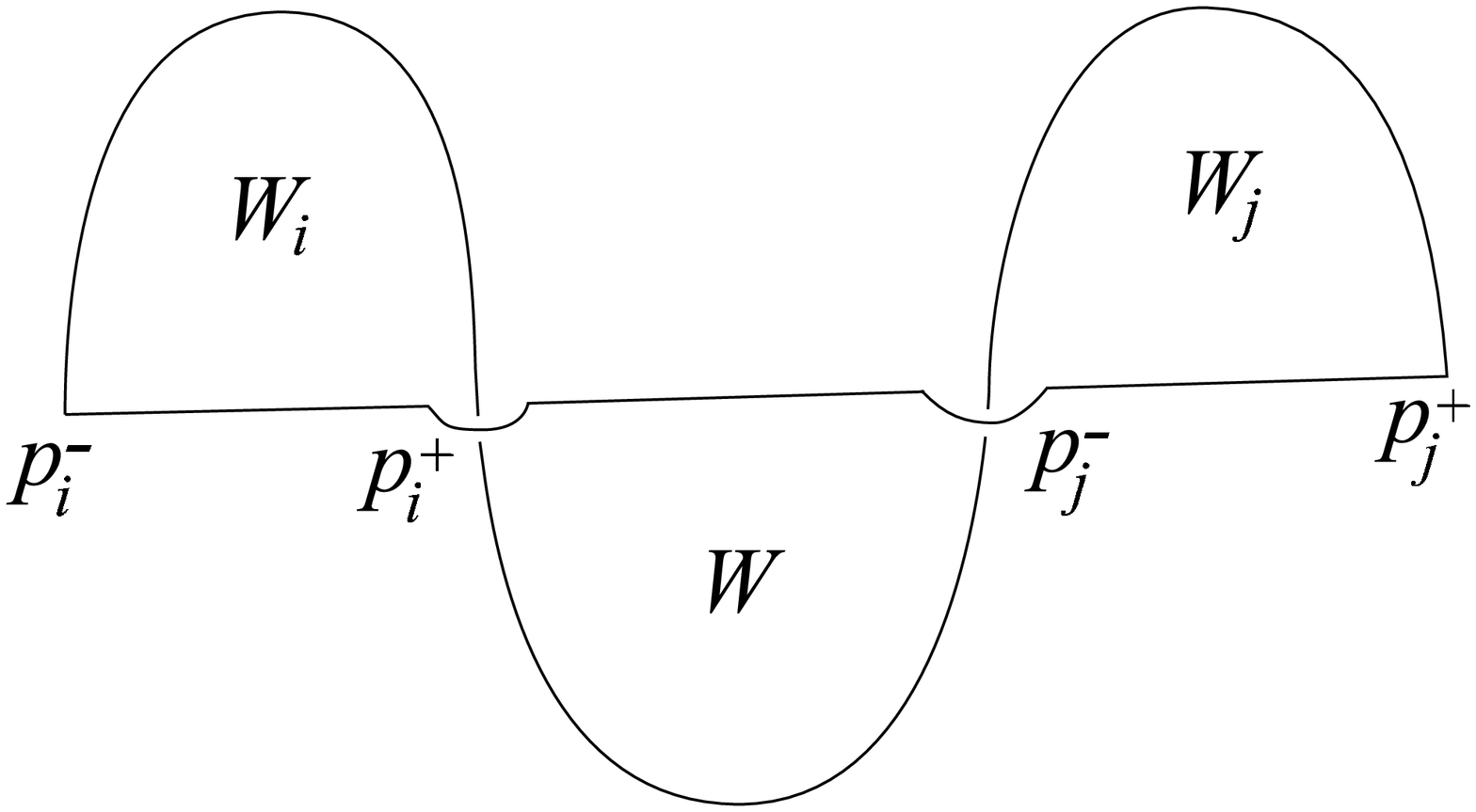}}
         \nocolon \caption{}\label{bandedWdisks}

\end{figure}

To show that $\tau(f)$ is well-defined it remains to check that
it does not depend on the choice of pairings of double points. If
$(p^{+}_{i},p^{-}_{i})$ and $(p^{+}_{j},p^{-}_{j})$
are paired by Whitney disks
$W_i$ and $W_j$ with $g_i=g_j$ then $(p^{+}_{i},p^{-}_{j})$ and
$(p^{+}_{j},p^{-}_{i})$ are also canceling pairs. Let $W$ be a
Whitney disk for $(p^{+}_{i},p^{-}_{j})$.
    A framed Whitney disk $W'$ for $(p^{+}_{j},p^{-}_{i})$ can be
formed by connecting $W$ to $W_i$ and $W_j$ using twisted strips
as in Figure~\ref{bandedWdisks} so that $$
I(W)+I(W')=I(W)+I(W_i)+I(W_j)-I(W)=I(W_i)+I(W_j). $$ Since any
two choices of pairings are related by a sequence of
interchanging pairs of double points in this way, it follows that
$\tau(f)$ does not depend on how the pairs are chosen from the
double points with the same group elements and opposite signs.

There is one more subtlety to check regarding the pairings which
is discussed in \cite{S} but neglected in \cite{FQ}: the pairing
of the {\it pre-images} of a canceling pair $(p^{+}_i,p^{-}_i)$
of double points with group element $g_i$ such that ${g_i}^2=1$.
Since $g_i={g_i}^{-1}$, the inverse image of the positive arc of
a Whitney disk $W_i$ can join an inverse image of $p^{-}_i$ to
\emph{either} of the two inverse images of $p^{+}_i$.

\begin{figure}[ht!] 
         \centerline{\includegraphics[width=0.4\hsize]{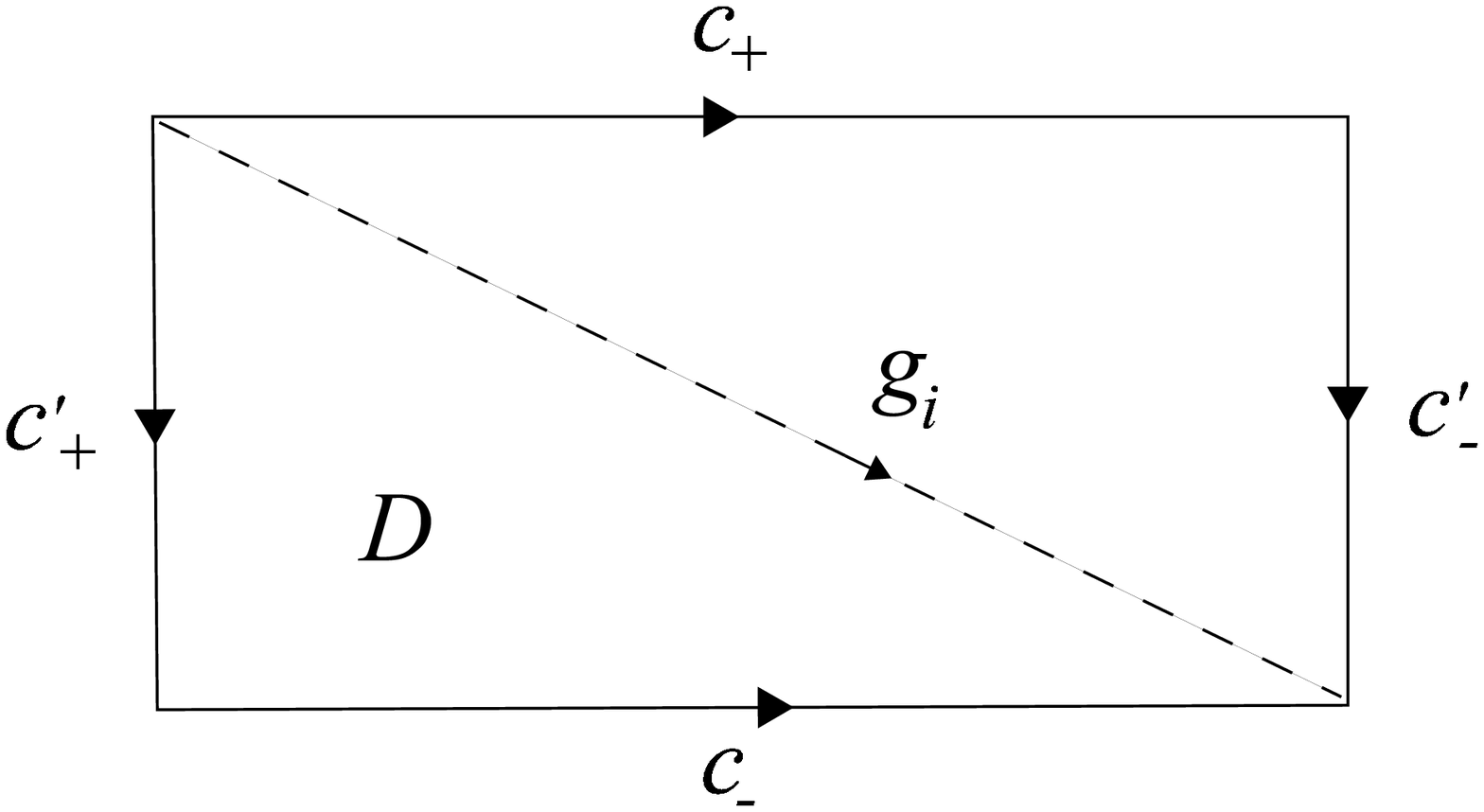}}
         \caption{The inverse image of the boundaries of two Whitneydisks
         for a canceling pair of doublepoints with group element $g_i$
         where $g_i^{2}=1$.}\label{RP2}

\end{figure}

Let $W_i$ and ${W_i}'$ be Whitney disks corresponding to the two
ways of pairing the inverse images of such a cancelling pair
$(p^{+}_{i},p^{-}_{i})$ with $g_i^2=1$. The union of the inverse
images of the boundary arcs of $W_i$ and ${W_i}'$ is a loop $c$
in $S^2$ which is the union of two pairs of arcs
$c_{\pm}:=f^{-1}(\partial_{\pm}W_i$) and
${c'}_{\pm}:=f^{-1}(\partial_{\pm}{W_i}')$
(see Figure~\ref{RP2}). By previous arguments
we may assume that $c$ is embedded
and bounds a 2-cell $D$ in $S^2$ such that $f$
restricts to an embedding on $D$. The union $A$ of the image of
$D$ together with $W_i$ and ${W_i}'$ is (after rounding corners)
an immersed $\R\PP^2$ representing $g_i$.  Since $W_i$ and
${W_i}'$ are correctly framed, the number of new intersections
between $A$ and $f$ that are created by perturbing $A$ to be
transverse to $f$ will be congruent to $\omega_2(A)$ modulo 2.
Each of these new intersections will have group elements $g_i$ or
1 so that $$ I(W_i)-I({W_i}')=(g_i,\lambda(f,A)+\omega_2(A)). $$
Thus $\tau(f)$ does not depend on the choice of the pairings of
the pre-images of the double points by the INT relation.

\begin{remn}\label{rem:2-torsion} If $\pi_1X$ has no 2-torsion then the above
immersion $\R\PP^2\imra X$ is spherical and hence the INT
relation only consists of intersections with spheres.
\end{remn}

We have shown that $\tau(f)$ is well-defined; it remains to show
that it is a homotopy invariant. As explained in
Section~\ref{sec:prelim} it suffices to show that it is invariant
ambient isotopies, finger moves, and (embedded) Whitney moves so
we will check that these moves do not change $\tau(f)$.  Any
isotopy of $f$ can be extended to the Whitney disks without
creating any new intersections between $f$ and the interiors of
the Whitney disks so that $\tau(f)$ is unchanged. A finger move
creates a cancelling pair of double points of $f$ equipped with a
\emph{clean} Whitney disk $W$, i.e. $W$ is embedded and $\int W
\cap f=\emptyset$. Since a finger move is supported in a
neighborhood of an arc it can be assumed to miss all pre-existing
Whitney disks. Thus $\tau(f)$ is unchanged by finger moves. A
Whitney move on $f$ pre-supposes the existence of a clean Whitney disk $W$.
We may assume that $W$ is included in any collection of
Whitney disks used to compute $\tau(f)$. The boundaries of all
other Whitney disks can be made disjoint from $\partial W$
by applying the move of
Figure~\ref{y-push} which does not change $\tau(f)$.
A Whitney move on $W$
eliminates the double points paired by $W$ and creates a pair of
new intersections between $f$ and $\int W_i$ for each point of
intersection in $\int W \cap \int W_i$. These new pairs of
intersections have cancelling contributions to $\tau(f)$ and so
the net change is zero. $\hfill \square$

\section{Proof of Theorem~\ref{thm2}}\label{sec:proof of thm2}

The ``if'' directions of Theorem~\ref{thm2} are clear from the
definition of $\tau(f)$. The ``only if'' direction will be shown
using the following lemma.

\begin{lem} \label{lem:I=0}
If $\tau(f)=0$ then after a homotopy of $f$ (consisting of finger
moves) the self-intersections of $f$ can be paired up by framed
Whitney disks $W_i$ with disjointly embedded boundaries such that
$I(W_i)=0 \in \Z [\pi_1 X \times \pi_1 X]$ for all $i$.
\end{lem}

\begin{figure}[ht!] 
         \centerline{\includegraphics[width=0.5\hsize]{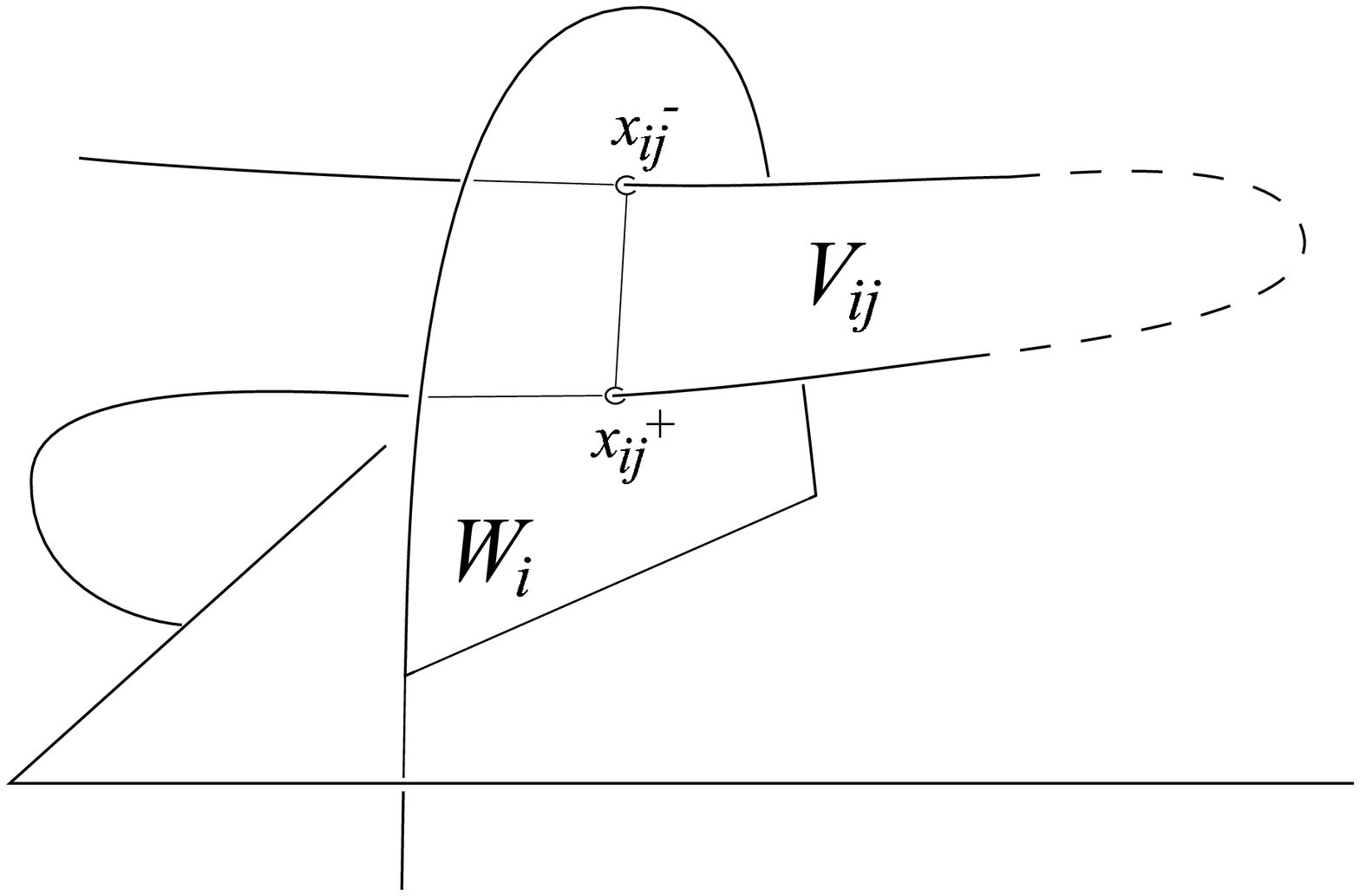}}
         \caption{A secondary Whitneydisk $V_{ij}$.}\label{secondaryWdisk}

\end{figure}

The geometric content of this lemma is that all the intersections
between $f$ and the interior of each Whitney disk $W_i$ are
paired by a second layer of Whitney disks: Since $I(W_i)=0$ the
intersections between $\int W_i$ and $f$ come in pairs
${x_{ij}}^\pm$ where $h_{{x_{ij}}^+}=h_{{x_{ij}}^-}\in \pi_1X$
and $\sign{{x_{ij}}^+}=-\sign{{x_{ij}}^-}$. The union of an arc
in $W_i$ (missing all double points of $W_i$) joining
${x_{ij}}^\pm$ and an arc in $f$ joining ${x_{ij}}^\pm$ (and
missing all double points of $f$) is a nullhomotopic loop which
bounds a Whitney disk $V_{ij}$ for the pair ${x_{ij}}^\pm$
(See Figure~\ref{secondaryWdisk}).

The proof of Lemma~\ref{lem:I=0} will be given shortly, but first
we use it to complete the proof of Theorem~\ref{thm2}.

We may assume, as just noted, that the self-intersections of $f$
are paired by framed Whitney disks $W_i$ with disjointly embedded
boundaries such that all intersections between the interiors of
the $W_i$ and $f$ are paired by Whitney disks $V_{ij}$. The
$V_{ij}$ can be assumed to be correctly framed after introducing
boundary twists (if necessary) around the arcs of the $\partial
V_{ij}$ that lie on the $W_i$.

The proof of Theorem~\ref{thm2} can be completed in two steps:
First use finger moves on $f$ to trade all intersections between
$f$ and the interiors of the $V_{ij}$ for new self-intersections
of $f$. These new self-intersections come paired by clean Whitney disks
disjoint from all other $W_i$. Next use the $V_{ij}$ to
guide Whitney moves on the $W_i$ eliminating all intersections
between $f$ and the interiors of the $W_i$. This second step may
introduce new interior intersections between Whitney disks but
these are allowed. These modified $W_i$ together with the new
clean Whitney disks have interiors disjoint from $f$ and
disjointly embedded boundaries.$\hfill \square$

$\hfill$

\begin{figure}[ht!] 
         \centerline{\includegraphics[scale=.43]{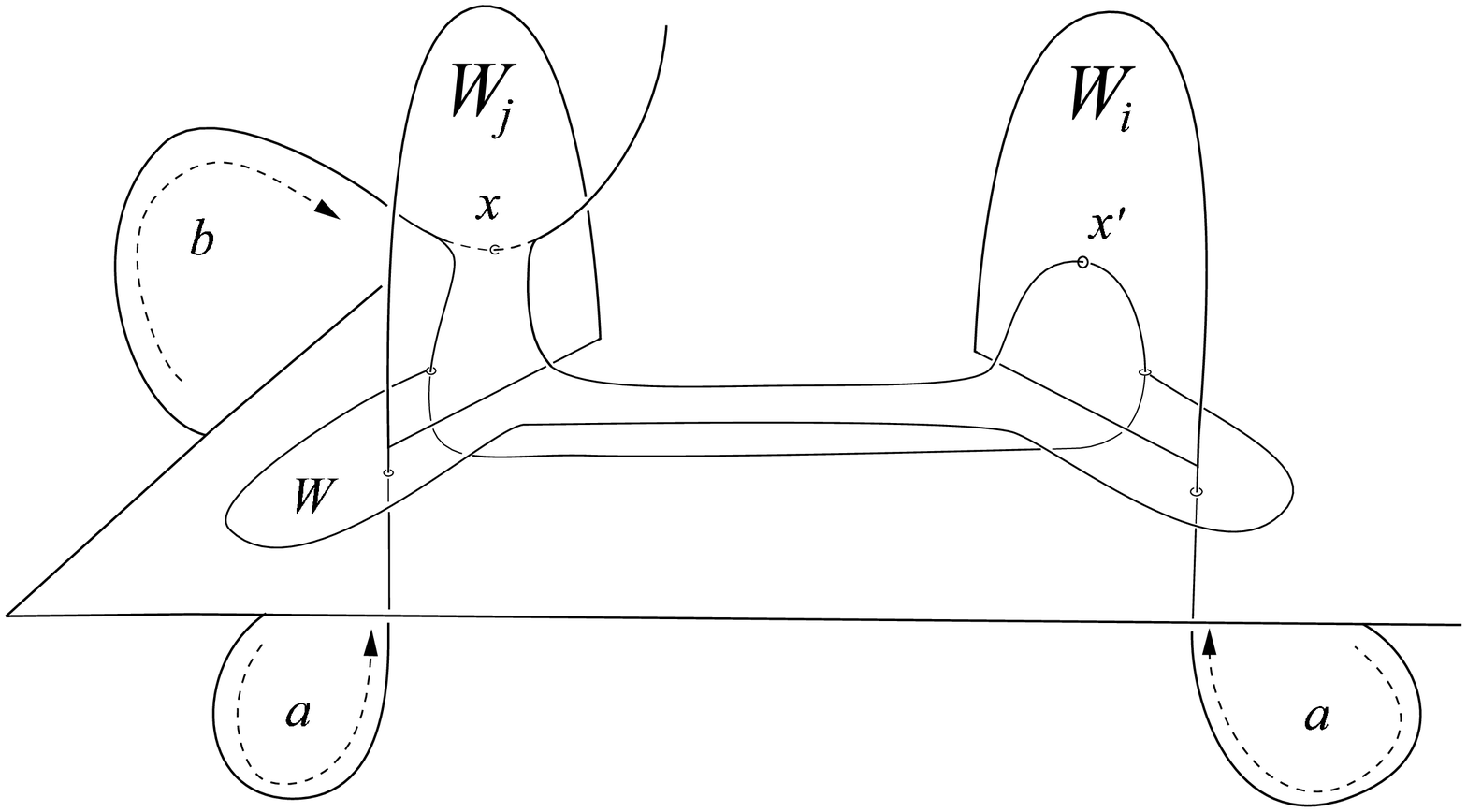}}
         \nocolon \caption{}\label{move}

\end{figure}

It remains to prove Lemma~\ref{lem:I=0}. The idea of the proof is
to first arrange for $\tau(f)$ to be given just in terms of
cancelling pairs of intersections between $f$ and the interiors
of the Whitney disks; then using the move described in Figure~\ref{move}
each cancelling pair can be arranged to occur on the same Whitney disk.

\begin{proof} Let $f$ satisfy $\tau(f)=0$ and $W_i$ be framed
Whitney disks pairing all the double points of $f$. The $W_i$ may
be assumed to have disjointly embedded boundaries after applying
the move of figure~\ref{y-push}. We now describe three modifications of $f$
and the collection of Whitney disks which can be used to
geometrically realize the relations FR, INT, and BC so that
$\tau(f)$ vanishes in the quotient of $\Z [\pi_1 X \times \pi_1 X]$
by the single
relation SC. (1) A finger move on $f$ guided by an arc
representing $a \in \pi_1(X)$ creates a cancelling pair of double
points of $f$ which are paired by a clean Whitney disk $W$. By
performing boundary twists and interior twists on $W$ one can
create intersections between $\int W$ and $f$ so that
$I(W)=n(a,1)+m(a,a)$ for any integers $n$ and $m$ such that
$n\equiv m$ modulo 2. (2) By similarly creating a clean Whitney disk $W$
and tubing into any immersed sphere representing $A\in
\pi_2(X)$ it can be arranged that
$I(W)=(a,\lambda(f,A)+\omega_2(A))$. (3) If a Whitney disk $W$
has an interior intersection point $x$ with $f$ that contributes
$\pm(a,b)$ to $I(W)$ then $x$ may be eliminated by a finger move
at the cost of creating a new pair $p^{\pm}$ of double points of
$f$ which admit a Whitney disk $W'$ (with embedded boundary
disjoint from existing Whitney disks) such that $\int {W'}$ has a
single intersection with $f$ and $I(W')=\mp(b,a)$ (See
Figure~\ref{pushab}). By using these three modifications we may assume that
our collection of Whitney disks satisfies $\sum_i I(W_i)=0$ in
$\Z [\pi_1 X \times \pi_1 X]$ modulo the SC relation.

\begin{figure}[ht!] 
         \centerline{\includegraphics[scale=.30]{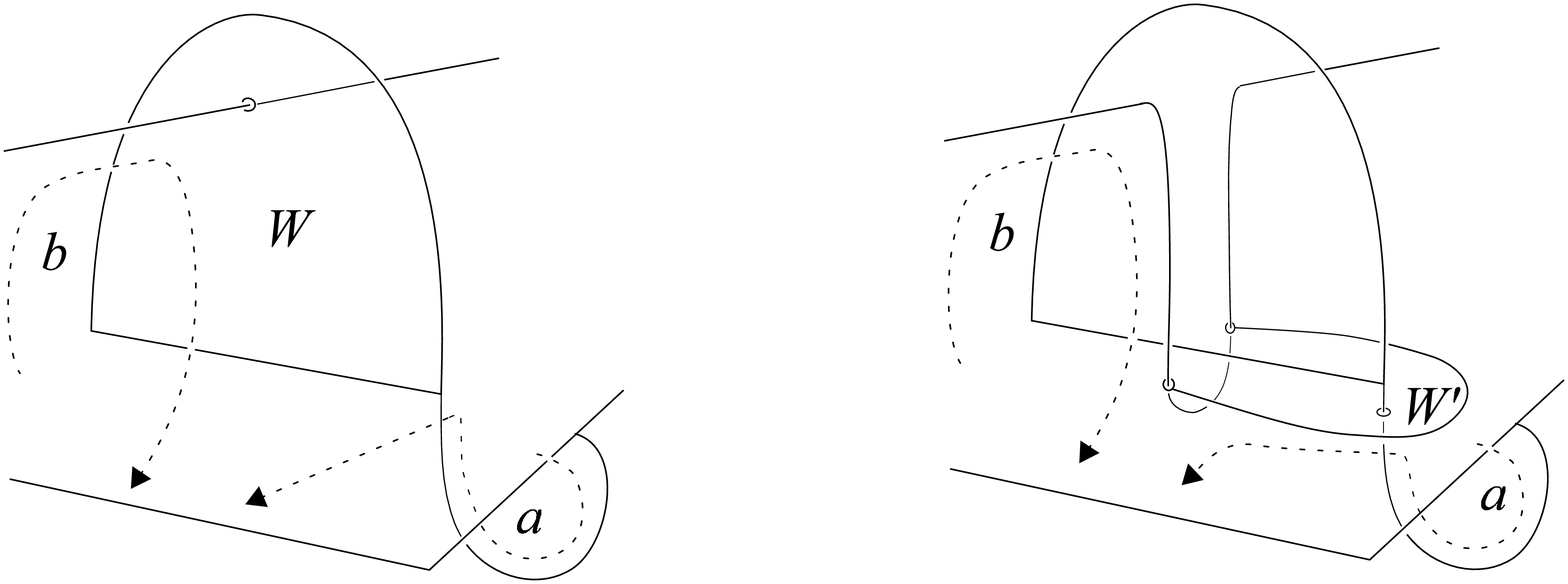}}
         \nocolon \caption{}\label{pushab}

\end{figure}

We can now move pairs of intersection points which have
algebraically cancelling contributions to $\tau(f)$ on to the
same Whitney disk as follows (see \cite{Y} for a detailed
description of the simply-connected case.) The finger move
illustrated in Figure~\ref{move} exchanges a point $x\in \int{W_j}\cap
f$ that contributes $(a,b)$ to $I(W_j)$ for a point $x'\in
\int{W_i}\cap f$ that contributes $(a,b)$ to $I(W_i)$. This
finger move also creates two new double points of $f$ which admit
a Whitney disk $W$ (with $\partial W$ embedded and disjoint from
all other Whitney disks) such that $I(W)=(b,a)-(b,a)=0$. By
performing this finger move through the negative arc of $W_j$
instead of the positive arc one can similarly exchange a point
$x\in \int{W_j}\cap f$ that contributes $-(a^{-1},ba^{-1})$ to
$I(W_j)$ for a point $x'\in \int{W_i}\cap f$ that contributes
$(a,b)$ to $I(W_i)$. In this way it can be arranged that all
double points of $f$ are paired by Whitney disks $W_i$ such that
$I(W_i)=0\in \Z [\pi_1 X \times \pi_1 X]$ for all $i$.

\end{proof}

\section{An invariant for a triple of immersed
spheres}\label{sec:triples}

In this section we define the cubic invariant $\lambda(f_1,f_2,f_3)$ of
Theorem~\ref{triples} and sketch the proof that it gives a
complete obstruction to making the $f_i$ disjoint. The invariant
is again given in terms of fundamental group elements determined
by secondary intersections, where in this case the relevant
intersections are between Whitney disks on two of the spheres and
the {\em other} sphere. However, we will pursue the point of view
of Remark~\ref{rem:Wall-type} and define the group elements via
Wall-type intersections between the $f_i$ and the Whitney disks.
While this approach initially increases the indeterminacy (due to
choosing whiskers for all the Whitney disks) it will eventually
serve to symmetrize the algebra and clarify the origin of the
$\s_3$-action in the invariant $\tau(f)$ for a single map of a
sphere. As before, we get an invariant taking values in a
quotient of $\Z[\pi \times \pi]$, this time via the
identification with $\Z[\pi \times \pi \times \pi / \Delta(\pi)]$
where $\Delta$ denotes the diagonal right action of
$\pi:=\pi_1X$.

Let $f_1,f_2,f_3:S^2\looparrowright X$ be an ordered triple of
oriented immersed spheres with pairwise vanishing Wall
intersections $\lambda(f_i,f_j)=0$ in an oriented 4-manifold $X$.
Choose Whitney disks with disjointly embedded boundaries pairing
all intersections between $f_i$ and $f_j$ for each pair $i\neq
j$. The notation for Wall intersections tacitly assumes that each
$f_i$ is equipped with a whisker (an arc connecting a basepoint
on $f_i$ to the basepoint of $X$). Now choose whiskers for each
of the Whitney disks. Orient all the Whitney disks as follows: If
$W^{ij}$ is a Whitney disk for a cancelling pair of intersections
between $f_i$ and $f_j$ with $i<j$ then take the positive (resp.
negative) arc of $W^{ij}$ to lie on $f_i$ (resp.\ $f_j$). As in
Section~\ref{sec:tau}, orient $W^{ij}$ by orienting $\partial
W^{ij}$ in the direction of the positive intersection point along
the positive arc then back to the negative intersection point
along the negative arc and taking a second outward-pointing
vector. To each intersection point $x$ between $f_k$ and the
interior of a Whitney disk $W^{ij}$ for a cancelling pair in $f_i
\cap f_j$ we associate three fundamental group elements as
follows: The {\em positive} (resp.\ {\em negative}) group element
is determined by a loop along the positive (resp.\ negative)
sheet, then back along $W^{ij}$ (and the whisker on $W^{ij}$).
The {\em interior} group element is determined by a loop along
$f_k$ to $x$ and back along $W^{ij}$. The three group elements
are ordered by the induced ordering of $\{i,j,k\}$ on the sheets.
Thus each such $x$ determines an element in $\Lambda:=
(\pi\times\pi \times \pi)/\Delta(\pi)$ where the diagonal right
action $\Delta$ is divided out in order to remove the choice of
the whisker for the Whitney disks. Denoting the positive,
negative and interior elements for $x\in \int W^{ij}_r \cap f_k$
by $g^{+}_r$, $g^{-}_r$ and $h_x$ respectively, we now set up
notation to measure the intersections between the spheres and the
Whitney disks by defining three elements in the abelian group
$\Lambda:=
   \Z[(\pi\times\pi \times \pi)/\Delta(\pi)]$ as follows:

$$
  I_3(W^{12}_r):= \sum_{x\in W^{12}_r \cap f_3}
\sign(x)(g^{+}_r,g^{-}_r,h_x)\in \Lambda,
$$

$$
I_2(W^{13}_r):= \sum_{x\in W^{13}_r \cap f_2}
-\sign(x)(g^{+}_r,g^{-}_r,h_x)\in \Lambda,
$$
and
$$
I_1(W^{23}_r):= \sum_{x\in W^{23}_r \cap f_1}
\sign(x)(g^{+}_r,g^{-}_r,h_x)\in \Lambda.
$$
  Denote by $\RR$ the
subgroup additively generated by $$ (a,b, \lambda(f_3,A)), (a,
\lambda(f_2,A),c), (\lambda(f_1,A),b,c)\in \Lambda $$ where
$a,b,c\in\pi$ and $A\in\pi_2X$ are arbitrary.

\begin{defi} \label{def:triple}
In the above setting define $$ \lambda(f_1,f_2,f_3):=\sum_r
I(W_r^{12})+\sum_r I(W_r^{23})+\sum_r I(W_r^{31})\in \Lambda/\RR.
$$ where the sums are over all Whitney disks for the
intersections between the $f_i$.
\end {defi}

\begin{remn}\label{triple example}
By modifying the construction of Section~\ref{sec:examples} one
can describe many triples with non-vanishing
$\lambda(f_1,f_2,f_3)$, for instance by shrinking three
components of the Bing double of the Hopf link in the complement
of the fourth component.
\end{remn}

Before sketching the proof of Theorem~\ref{triples} we now
describe a nice formalism which explains the presence of the
$\s_3$ indeterminacy in the definition of $\tau(f)$ which is
absent in the case of $\lambda(f_1,f_2,f_3)$ for a triple. In
both cases one assigns (two respectively three) fundamental group
elements to each intersection point between the interior of a
Whitney disk and a sheet of a sphere which we will refer to as
the {\em interior} sheet.

\begin{figure}[ht!] 
         \centerline{\includegraphics[width=.8\hsize]{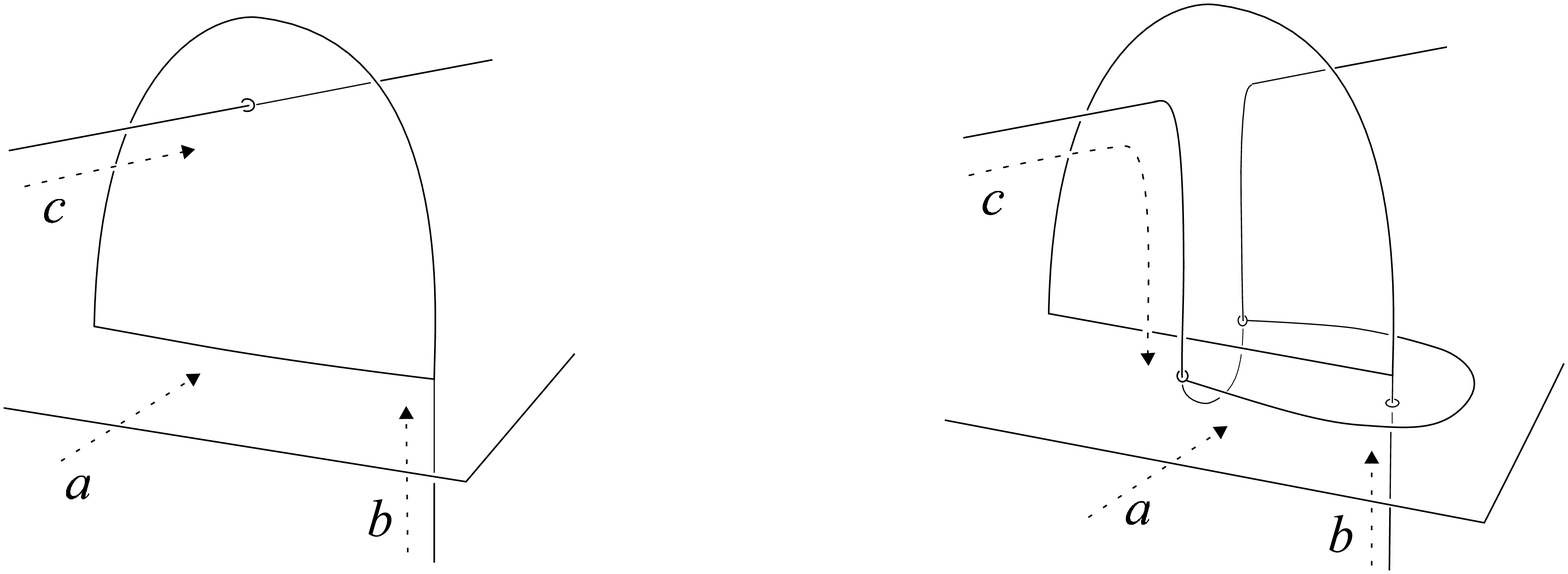}}
         \nocolon \caption{}\label{pushabc}

\end{figure}

\begin{figure}[ht!] 
         \centerline{\includegraphics[width=.8\hsize]{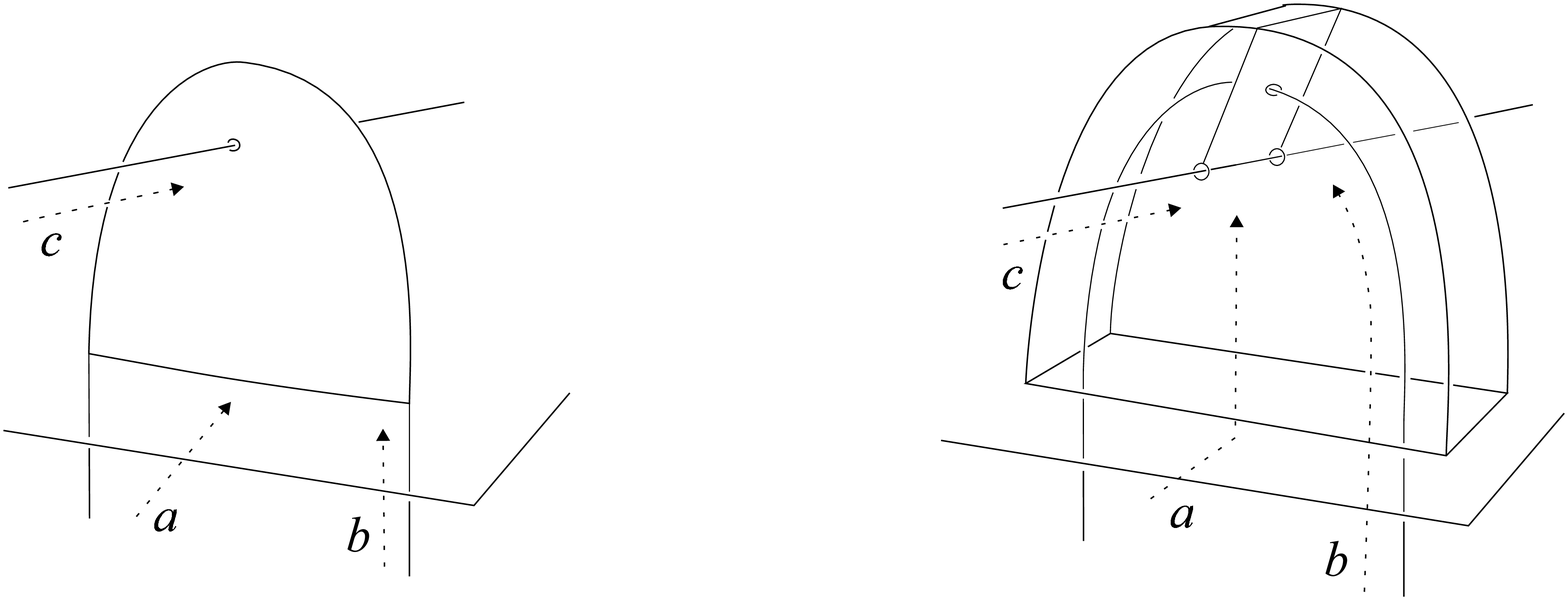}}
         \nocolon \caption{}\label{Whitneymove}

\end{figure}

For each such intersection point the corresponding interior sheet
``interacts'' with the positive and negative sheets of the
Whitney disk in the following sense: By pushing down the interior
sheet into the positive (resp.\ negative) sheet one can eliminate
the original intersection point at the cost of creating a new
cancelling pair of intersections which admits a new Whitney disk
which has an interior intersection point with the negative (resp.
positive) sheet (see Figure~\ref{pushabc}). Note that this trading of one
intersection point for another takes place in a neighborhood of
the original Whitney disk and the effect of pushing down into a
sheet is the same as the effect of doing a Whitney move (see
Figure~\ref{Whitneymove}). It is clear that any invariant defined in terms of
such intersections will have ``local'' indeterminacies
corresponding to this local interaction between the three sheets.

\begin{figure}[ht!] 
         \centerline{\includegraphics[scale=.30]{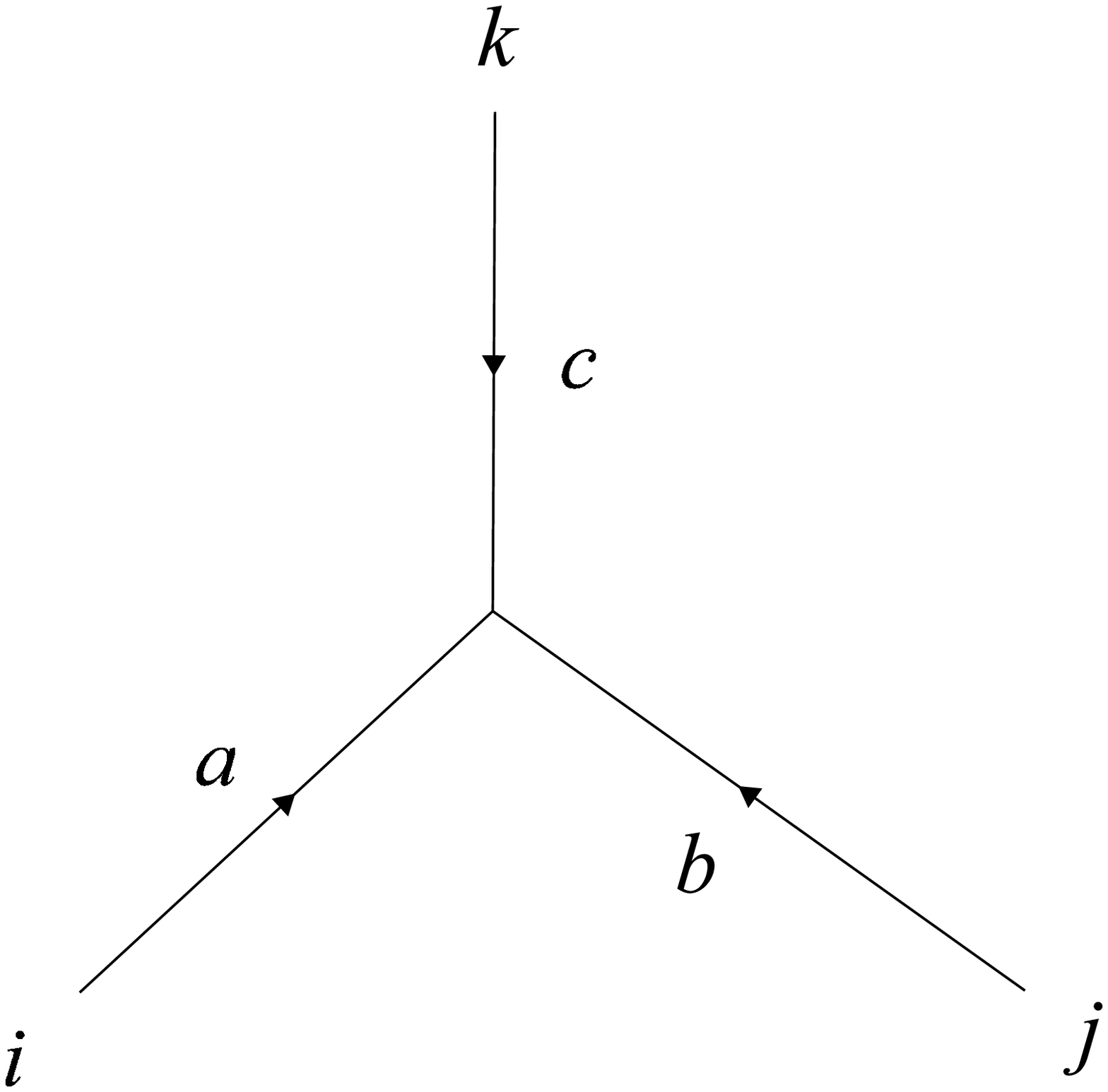}}
         \nocolon \caption{}\label{tree}

\end{figure}

This interaction can be described by associating a decorated
uni-trivalent tree with one interior vertex to each intersection
point $x$ between the interior of a Whitney disk and a sheet as
follows (see Figure~\ref{tree}). The interior vertex represents the
Whitney disk and the three univalent vertices represent sheets of
$f_i$, $f_j$ and $f_k$, two of which are the positive and
negative sheets of the Whitney disk, the other being the interior
sheet corresponding to $x$.
   The three
edges
   are oriented inward and represent the corresponding
positive, negative and interior group elements.
    {\em The relations that are
forced on the triple of group elements by the above described
interactions between the sheets correspond to (signed) graph
automorphisms of the tree which preserve the labels of the
univalent vertices (with the sign given by the sign of the
induced permutation of the vertices)}.

   In particular,
if $i=j=k$ as in the case of $\tau(f)$, then the automorphism
group is $\s_3$. In fact, if one defines $\tau(f)$ by the
conventions of this section (i.e., choosing whiskers for the
Whitney disks, etc.) then the $\s_3$-action is obvious because
the SC and BC relations become $(a,b,c)=-(b,a,c)$ and
$(a,b,c)=-(a,c,b)$ respectively. If $i$, $j$ and $k$ are distinct
as in the case of $\lambda(f_1,f_2,f_3)$ then the automorphism
group is trivial and no ``local'' relations are needed. We will
see in Section~\ref{sec:general tau} that the expected $\s_2$
indeterminacy (corresponding to graph automorphisms which fix one
univalent vertex) is present in the definition of $\tau(f_1,f_2)$
for a pair.

\begin{remn}\label{Jacobi}
Higher order generalizations of $\tau$ will have indeterminacies
that correspond to antisymmetry and Jacobi relations known from
the theory of finite type invariants of links in $3$-manifolds.
\end{remn}

We now give a proof of Theorem~\ref{triples}. Since the arguments are
essentially the same as those of Section~\ref{sec:tau
well-defined} (with the added convenience of working with much
less indeterminacy) the steps will be indicated with details
omitted.

\begin{proof}
First note that since the $f_i$ are ordered all the Whitney disks
are canonically oriented via our convention. Thus the signs
associated to the intersections between the Whitney disks and the
$f_i$ are well-defined. Also, the element of $\Lambda$ associated
to such an intersection point does not depend on the choice of
whisker for the Whitney disk because we are modding out by the
diagonal right action $\Delta(\pi)$ of $\pi$. Since any two
Whitney disks with the same boundary differ by an element $A \in
\pi_2 X$, $\lambda(f_1,f_2,f_3)$ does not depend on the choices
of the interiors of the Whitney disks because it is measured in
the quotient of $\Lambda$ by $\RR$.

In order to show that $\lambda(f_1,f_2,f_3)$ does not depend on
the boundaries of the Whitney disks it is convenient to
generalize the definition of $\lambda(f_1,f_2,f_3)$ to allow weak
Whitney disks along the lines of Section~\ref{sec:tau}. For each
$y\in \partial W^{ij}\cap \partial W^{ik}$, where $j\neq k$ and
$(\overrightarrow{\partial W^{ij}},\overrightarrow{\partial
W^{ik}})$ agrees with the orientation of $f_i$ at $y$, define
$J(y)\in \Lambda$ from the following three group elements: Give
$W^{ij}$ and $W^{ik}$ a common whisker at $y$ then take the
positive and negative group elements of $W^{ij}$ together with
the group element of $W^{ik}$ corresponding to the sheet $f_k$.
Define the sign of $J(y)$ to be equal to the sign of the
permutation $(ijk)$.
   The generalized definition of $\lambda(f_1,f_2,f_3)$
includes the sum of $J(y)$ over such $y$. Note that this
generalization does not require any new relations and reduces to
the original definition after eliminating the intersections and
self-intersections between boundaries of the Whitney disks in the
usual way (Figure~\ref{y-push}). The arguments of Section~\ref{sec:tau
well-defined} now apply to show that $\lambda(f_1,f_2,f_3)$ does
not depend on the boundaries of the Whitney disks: pushing the
boundary of a Whitney disk across an intersection point
(Figure~\ref{pushacross}) creates
an interior intersection in the collar of the Whitney disk and a
boundary intersection with cancelling contributions to
$\lambda(f_1,f_2,f_3)$. Independence of the choice of pairings of
the intersection points (see Figure~\ref{bandedWdisks}) also follows.
Note that in
the present setting there are no subtleties concerning the
pre-images of intersection points. Thus $\lambda(f_1,f_2,f_3)$ is
well-defined.

To see that $\lambda(f_1,f_2,f_3)$ is invariant under homotopies
of the $f_i$ it suffices to check invariance under finger moves
and Whitney moves on embedded Whitney disks. Finger moves only
create embedded Whitney disks which clearly don't change
$\lambda(f_1,f_2,f_3)$. A Whitney move on an embedded Whitney disk
$W$ which eliminates a pair of self-intersections on $f_i$
will create a new pair of intersections between $f_i$ and
$W^{jk}$ for each intersection between $W^{jk}$ and $W$, but
these new intersections have cancelling contributions to
$\lambda(f_1,f_2,f_3)$ which remains unchanged. The same applies
to a Whitney move that eliminates a pair of intersections between
two different spheres and since the embedded Whitney disk can be
assumed to have been included in any collection used to compute
$\lambda(f_1,f_2,f_3)$ again the invariant is unchanged. Thus
$\lambda(f_1,f_2,f_3)$ only depends on the homotopy classes of
the $f_i$.

It follows directly from the definition that $\lambda(f_1,f_2,f_3)$
satisfies properties (i), (ii), (iv) and (v)
of Theorem~\ref{triples}.
Property (iii) can be checked as follows:
$\lambda(f,f,f)=\lambda(f_1,f_2,f_3)$ where the $f_i$ are
parallel copies of $f$. Since the normal bundle
of $f$ is trivial, each self-intersection of $f$ gives nine
intersection points among the $f_i$ of which exactly three are self-intersections.
Thus each Whitney disk $W$ for a canceling pair of
{\em self}-intersections of $f$ yields six canceling pairs
of intersections among the $f_i$
paired by Whitney disks which are essentially
parallel copies of $W$.
If $W$ has an interior intersection with $f$ that contributes
$(a,b,c)$ to $\tau(f)$ (expressed in the notation of this section)
then there will be six corresponding terms contributing
$\sum_{\sigma\in \s_3}(a,b,c)^\sigma$ to
$\lambda(f_1,f_2,f_3)$. See Figure~\ref{parallelsheets} for a schematic
illustration where we have circled the points corresponding to
one of the contributions, an intersection between $f_3$ and a Whitney disk
on $f_1$ and $f_2$.

\begin{figure}[ht!] 
         \centerline{\includegraphics[width=0.4\hsize]{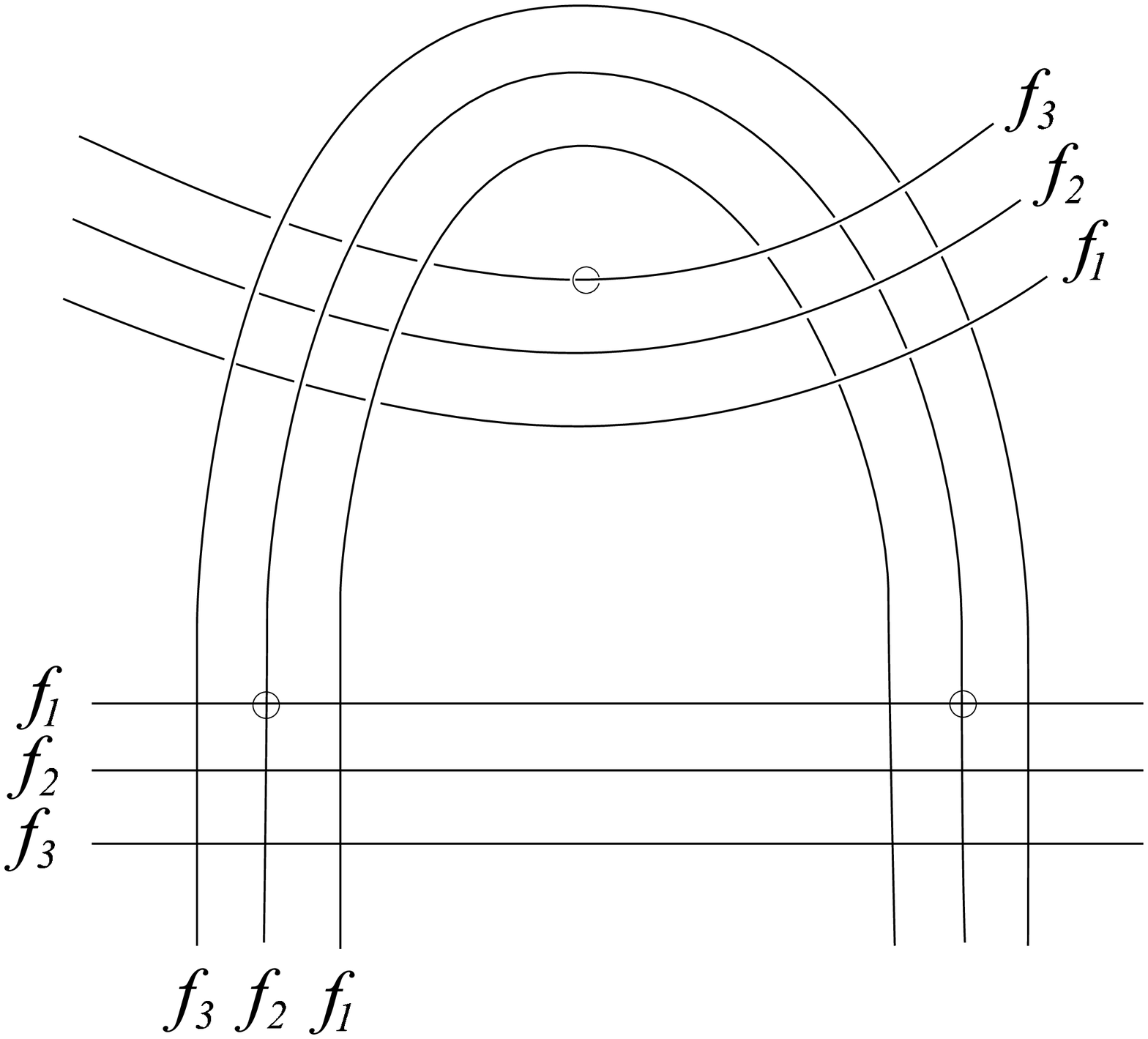}}
         \nocolon \caption{}\label{parallelsheets}

\end{figure}

To complete the proof of  Theorem~\ref{triples} we now show that
the $f_i$ can be homotoped to pairwise disjoint maps if
$\lambda(f_1,f_2,f_3)=0$ (the converse is clear). First use
$\lambda(f_1,f_2)=0$ to separate $f_1$ and $f_2$ by pushing $f_2$
off the Whitney disks $W^{12}$ and then doing Whitney moves on
$f_1$. Now $\lambda(f_1,f_2,f_3)=0$ is given completely in terms
of Whitney disks pairing $f_3 \cap f_1$ and Whitney disks pairing
$f_3 \cap f_2$. The arguments of Lemma~\ref{lem:I=0} (tubing the
Whitney disks into spheres and using the move of Figure~\ref{move}) can
now be applied to get a second layer of Whitney disks which pair
all intersections between $f_1$ and the Whitney disks for $f_3
\cap f_2$. After pushing any intersections between the secondary
Whitney disks and $f_1$ down into $f_1$, the secondary Whitney disks
can be used to make $f_1$ disjoint from the Whitney disks
for $f_3 \cap f_2$. After pushing down any intersections with
$f_2$, the Whitney disks on $f_3 \cap f_2$ may now be used to
eliminate all intersections between $f_2$ and $f_3$. After
similarly applying the arguments of Lemma~\ref{lem:I=0} to get
secondary Whitney disks pairing the intersections between $f_2$
and the Whitney disks for $f_3 \cap f_1$, one can push down
intersections and do Whitney moves to eliminate all intersections
between $f_1$ and $f_3$ so that the $f_i$ are pairwise disjoint.
\end{proof}

\section{An invariant for $n$ immersed
spheres}\label{sec:general tau}

In this section we define the invariant $\tau(f_1, \ldots, f_n)$
of Theorem~\ref{n spheres} which obstructs homotoping $n$
immersed spheres in a 4-manifold $X$ to disjoint embeddings and
is the complete obstruction to raising the height of a
Whitney-tower. As before, the invariant is determined by
intersections between the spheres and their Whitney disks. We
continue along the lines of Section~\ref{sec:triples} by working
with {\em based} Whitney disks, i.e. choosing a whisker for each
Whitney disk, and identifying $\pi\times\pi$ with the quotient
$(\pi\times\pi \times \pi)/\Delta(\pi)$, where $\Delta$ denotes
the diagonal {\em right} action of $\pi:=\pi_1X$. Since there are
now $n$ different choices for each of the three sheets
interacting at any Whitney disk, $\tau(f_1,\ldots,f_n)$ will take
values in ${n \choose 1}+2{n \choose 2}+{n \choose 3}$ copies of
$\Z[(\pi\times\pi \times \pi)/\Delta(\pi)]$ modulo some relations
which are generalized versions of the SC, BC, FR and INT
relations of Section~\ref{sec:tau}. The ideas of this section are
completely analogous to those of previous sections, the only
novelty being the need to develop notation and conventions to
handle $n$ maps. Proofs will be omitted.

Let $f_1, \dots, f_n: S^2\looparrowright X^4$ be a collection of
oriented immersed spheres with vanishing primary $\mu$ and
$\lambda$ intersection numbers. Choose based, framed Whitney disks
with disjointly embedded boundaries pairing all
intersections and self-intersections among the $f_i$. Orient all
the Whitney disks as follows: If $W^{ij}$ is a Whitney disk for a
cancelling pair of intersections between $f_i$ and $f_j$ with
$i\leq j$ then take the positive arc of $W^{ij}$ to lie on $f_i$.
Orient $W^{ij}$ by orienting $\partial W^{ij}$ in the direction
of the positive intersection point along the positive arc and
taking a second outward-pointing vector. To each intersection
point $x$ between $f_k$ and the interior of a Whitney disk
$W^{ij}$ for a cancelling pair in $f_i \cap f_j$ we associate
three fundamental group elements corresponding to the three
sheets as follows: The {\em positive} (resp.\ {\em negative})
group element is determined by a loop along the positive (resp.
negative) sheet, then back along $W^{ij}$ (and the whisker on
$W^{ij}$). The {\em interior} group element is determined by a
loop along $f_k$ to $x$ and back along $W^{ij}$. The three group
elements are ordered by the induced ordering of the maps on the
sheets together with the convention that the positive element
precedes the negative element which precedes the interior
element. Thus each such $x$ determines an element in $
\Z[(\pi\times\pi \times \pi)/\Delta(\pi)]$ where the sign is
given by $\sign (x)$ times the sign of the permutation $(ijk)$
with the above ordering conventions.

Denoting the positive, negative and interior elements for $x\in
\int W^{ij}_r \cap f_k$ by $g^{+}_r$, $g^{-}_r$ and $h_x$
respectively, we now set up notation to measure the intersections
between the spheres and the Whitney disks. For $i<j<k$ define $$
I_i(W^{ii}_r):= \sum_{x\in W^{ii}_r \cap f_i}
\sign(x)(g^{+}_r,g^{-}_r,h_x)_{iii}\in \Lambda_{iii} $$ $$
I_j(W^{ii}_r):= \sum_{x\in W^{ii}_r \cap f_j}
\sign(x)(g^{+}_r,g^{-}_r,h_x)_{iij}\in \Lambda_{iij} $$ $$
I_i(W^{ij}_r):= \sum_{x\in W^{ij}_r \cap f_i}
(-1)\sign(x)(g^{+}_r,h_x,g^{-}_r)_{iij}\in \Lambda_{iij} $$ $$
I_j(W^{ij}_r):= \sum_{x\in W^{ij}_r \cap f_j}
\sign(x)(g^{+}_r,g^{-}_r,h_x)_{ijj}\in \Lambda_{ijj} $$ $$
I_i(W^{jj}_r):= \sum_{x\in W^{jj}_r \cap f_i}
\sign(x)(h_x,g^{+}_r,g^{-}_r)_{ijj}\in \Lambda_{ijj} $$ $$
I_k(W^{ij}_r):= \sum_{x\in W^{ij}_r \cap f_k}
\sign(x)(g^{+}_r,g^{-}_r,h_x)_{ijk}\in \Lambda_{ijk} $$ $$
I_j(W^{ik}_r):= \sum_{x\in W^{ik}_r \cap f_j}
(-1)\sign(x)(g^{+}_r,h_x,g^{-}_r)_{ijk}\in \Lambda_{ijk} $$ $$
I_i(W^{jk}_r):= \sum_{x\in W^{jk}_r \cap f_i}
\sign(x)(h_x,g^{+}_r,g^{-}_r)_{ijk}\in \Lambda_{ijk}. $$ Denote
by $\Lambda$ the direct sum $$ \bigoplus_i \Lambda_{iii}
\bigoplus_{i<j}(\Lambda_{iij} \oplus \Lambda_{ijj})
\bigoplus_{i<j<k} \Lambda_{ijk}. $$ where each $ \Lambda_{abc}$
is a copy of the abelian group $\Z[(\pi\times\pi \times
\pi)/\Delta(\pi)]$.
\begin{defi}
In the above setting define $$ \tau(f_1, \ldots, f_n):=\sum
I_i(W^{jk}_r)\in \Lambda/\RR $$ where the sum is over $i$ and all
the Whitney disks $W^{jk}_r$ ($j\leq k$). The relations $\RR$ are
additively generated by the following equations:

For all $i$ and $j$ with $i\leq j$,
\begin{eqnarray}
    (a,b,c)_{iij} & = & -(b,a,c)_{iij},\nonumber\\
    (a,b,c)_{ijj} & = & -(a,c,b)_{ijj} ,\nonumber\\
   (a,a,b)_{iij} &= & (a,b,b)_{ijj},\nonumber
\end{eqnarray}
\begin{eqnarray}
\sum_{k\leq i}(\lambda(f_k,A),a,b)_{kij}
-\sum_{i<k<j}(a,\lambda(f_k,A),b)_{ikj} +\sum_{k\geq
j}(a,b,\lambda(f_k,A))_{ijk}\kern-1.5in&&\nonumber\\
&& +(\omega_2A)(a,b,b)_{ijj} =0,\nonumber
\end{eqnarray}
where the sums are over $k$. Here $a,b \in \pi$ and $A\in \pi_2X$
are arbitrary. As in Section~\ref{sec:tau}, $A$ may be any
immersed $\R\PP^2$ representing $ab^{-1}$ whenever $i=j$ and
$ab^{-1}$ is of order two.
\end{defi}

\begin{remn}
The first three equations give local relations corresponding to
the interaction of the sheets at a Whitney disk (the third
equation corresponds to the boundary-twist operation). The fourth
equation gives global relations corresponding to indeterminacies
due to changing the homotopy class of a Whitney disk by tubing
into $A$.
\end{remn}

\begin{remn}
Note that $\tau(f_1,\dots,f_n)$ reduces to the invariant $\tau(f)$ of
Section~\ref{sec:tau} in the case $n=1$ (via the map
$(a,b,c)\longmapsto (ba^{-1},ca^{-1})$). Also, by ignoring all
terms with any non-distinct indices in the case $n=3$ we recover
the invariant $\lambda(f_1,f_2,f_3)$ of
Section~\ref{sec:triples}.
\end{remn}

\Addresses\recd


\begin{thebibliography}


\bibitem{COT} {\bf T Cochran}, {\bf K Orr}, {\bf P Teichner}, {\it Knot
concordance, Whitney towers and $L^2$-signatures}. Preprint
1999, posted at http://xxx.lanl.gov/abs/math/9908117

\bibitem{FQ}  {\bf M\,H Freedman}, {\bf F Quinn}, {\it The topology of
4-manifolds}, Princeton Math. Series 39, Princeton, NJ, 1990

\bibitem{GL} {\bf S Garoufalidis}, {\bf J Levine}, {\it Homology surgery
and invariants of 3-manifolds},
Preprint 2000, see http://xxx.lanl.gov/abs/math.GT/0005280

\bibitem{GR} {\bf S Garoufalidis}, {\bf L Rozansky},
{\it The loop expansion of the Kontsevich integral, abelian invariants of
knots and S-equivalence}.
Posted at http://xxx.lanl.gov/abs/math.GT/0003187

\bibitem{HK} {\bf I Hambleton}, {\bf M Kreck}, {\it Cancellation of
hyperbolic forms and topological four-manifolds},  J. Reine
Angew. Math. 443 (1993) 21--47

\bibitem{HS} {\bf W\,C Hsiang}, {\bf R\,H Szczarba}, {\it On embedding
surfaces in four-manifolds},
Proc. Sympos. Pure Math. 22 (1971) 97--103

\bibitem{KM} {\bf M Kervaire}, {\bf J Milnor}, {\it On 2-spheres in 4-manifolds},
Proc. Nat. Acad. Sci., Vol. 47, (1961) 1651--1657

\bibitem{LW} {\bf R Lee}, {\bf D Wilczynski}, {\it Locally flat $2$-spheres in simply
connected $4$-manifolds},  Comment. Math. Helv. 65, no. 3
(1990) 388--412

\bibitem{L} {\bf G\,S Li},  {\it An invariant of link homotopy in dimension
four }, Topology 36 (1997), 881-897

\bibitem{M} {\bf Y Matsumoto}, {\it Secondary intersectional properties of
4-manifolds and Whitney's trick.}, Proc. Symp. Pure Math. 32
Part II (1978) 99-107

\bibitem{R}  {\bf V Rohlin}, {\it Two-dimensional submanifolds of four dimension
manifolds }, J. Func. Anal. Appl. 5 (1971), 39-48


\bibitem{S} {\bf R Stong}, { \it Existence of $\pi_1$-negligible
embeddings in $4$-manifolds: A correction to Theorem 10.5 of
Freedman and Quinn }, Proc. of the A.M.S. 120 (4) (1994)
1309-1314

\bibitem{T} {\bf A\,G Tristam}, {\it Some coborbism invariants for links},
Proc. Cambridge Philos. Soc. 66 (1969) 251--264

\bibitem{W}  {\bf C\,T\,C\, Wall}, { \it Surgery on Compact Manifolds},
London Math.Soc.Monographs~1, Academic Press~1970 or Second Edition,
edited by A. Ranicki, Math. Surveys and Monographs 69, A.M.S

\bibitem{Wh} {\bf H Whitney}, {\it The self intersections of a smooth
$n$-manifold in $2n$-space}, Annals of Math. 45 (1944) 220-246

\bibitem{Y} {\bf M Yamasaki}, {\em Whitney's trick for three 2-dimensional
homology classes of 4-manifolds}, Proc. Amer. Math. Soc. 75
(1979) 365-371

\end{thebibliography}
\end{document}